\documentclass[aps,prl,twocolumn,superscriptaddress,amsmath,amsfonts,amssymb]{revtex4-1}

\usepackage{bm}
\usepackage{graphicx}
\usepackage{subfigure}
\usepackage{array}
\usepackage{amssymb}
\usepackage{amsfonts}
\usepackage{amsmath}
\usepackage{mathrsfs}
\usepackage{color}
\usepackage{booktabs}
\usepackage{threeparttable}
\usepackage{multirow}
\usepackage{epsfig}
\usepackage{threeparttable}
\usepackage{chngpage}
\usepackage{latexsym}
\usepackage[citecolor=blue,linkcolor=blue,colorlinks=true]{hyperref}%
\usepackage{epstopdf}
\usepackage{framed}
\definecolor{shadecolor}{rgb}{0.92,0.92,0.92}

\allowdisplaybreaks[4]

\begin{document}

\title{Continuity scaling: A rigorous framework for detecting and quantifying causality accurately}

\author{Xiong Ying}
\affiliation{School of Mathematical Sciences, SCMS, and SCAM, Fudan University, Shanghai 200433, China}
\affiliation{Research Institute for Intelligent Complex Systems,
CCSB, and LCNBI, Fudan University, Shanghai 200433, China}
\affiliation{State Key Laboratory of Medical Neurobiology, and MOE Frontiers Center for Brain Science, Institutes of Brain Science, Fudan University, Shanghai 200032, China}

\author{Si-Yang Leng}
\affiliation{Research Institute for Intelligent Complex Systems,
CCSB, and LCNBI, Fudan University, Shanghai 200433, China}
\affiliation{Institute of AI and Robotics, Academy for Engineering and Technology, Fudan University, Shanghai 200433, China}

\author{Huan-Fei Ma}
\affiliation{School of Mathematical Sciences, Soochow University, Suzhou 215006, China}

\author{Qing Nie}
\affiliation{Department of Mathematics, Department of Developmental and Cell Biology, and NSF-Simons Center for Multiscale Cell Fate Research, University of California, Irvine, CA 92697-3875, USA}

\author{Ying-Cheng Lai}
\affiliation{School of Electrical, Computer, and Energy Engineering, Arizona State University, Tempe, Arizona 85287-5706, USA}

\author{Wei Lin} \email{wlin@fudan.edu.cn}
\affiliation{School of Mathematical Sciences, SCMS, and SCAM, Fudan University, Shanghai 200433, China}
\affiliation{Research Institute for Intelligent Complex Systems,
CCSB, and LCNBI, Fudan University, Shanghai 200433, China}
\affiliation{State Key Laboratory of Medical Neurobiology, and MOE Frontiers Center for Brain Science, Institutes of Brain Science, Fudan University, Shanghai 200032, China}
\affiliation{Shanghai Artificial Intelligence Laboratory, Shanghai 200232, China}

\date{\today}

\begin{abstract}
Data based detection and quantification of causation in complex, nonlinear dynamical systems is of paramount importance to science, engineering and beyond.
Inspired by the widely used methodology in recent years, the cross-map-based techniques,
we develop a general framework to advance towards a comprehensive understanding of dynamical causal mechanisms,
which is consistent with the natural interpretation of causality. In particular, instead of measuring the smoothness of the cross map as conventionally implemented, we define causation through measuring the {\it scaling law} for the continuity of the investigated dynamical system directly.
The uncovered scaling law enables accurate, reliable, and efficient detection of causation and assessment of its strength in general complex dynamical systems, outperforming those existing representative methods.
The continuity scaling based framework is rigorously established and demonstrated using datasets from model complex systems and the real world.
\end{abstract}

\maketitle

\section{Introduction}

Identifying and ascertaining causal relations is a problem of paramount
importance to science and engineering with broad
applications~\cite{bunge2017causality,pearl2009causality,Runge2019inferring}. For example,
accurate detection of causation is key to identifying the origin of diseases
in precision medicine~\cite{CV:2015} and is important to fields such
as psychiatry~\cite{SSFRLPWBBA:2016}. Traditionally, associational concepts are often misinterpreted as causation~\cite{CoxHinkley1979,CoverThomas2012}, while causal analysis in fact goes one step further beyond association in a sense that, instead of using static conditions, causation is induced under changing conditions~\cite{JPearl2009}. The principle of Granger
causality formalizes a paradigmatic framework~\cite{Wiener:1956,Granger:1969,Haufe2013}, quantifying causality in terms of prediction improvements, but, because of its linear, multivariate, and statistical regression nature,
the various derived methods require extensive data~\cite{ding2006granger}. Entropy based methods~\cite{Schreiber:2000,frenzel2007partial,
vicente2011transfer,runge2012escaping,SCB:2014,CSB:2015,
STB:2015,research2021} face a similar difficulty. Another issue with the Granger causality
is the fundamental requirement of separability of the underlying dynamical
variables, which usually cannot be met in the real world systems. To overcome these
difficulties, the cross-map-based techniques, paradigms tailored to dynamical systems, have been developed and have gained widespread attention in the past decade~\cite{PRE_Hirata2010,PRE_Quiroga2000,PhyD_Arnhold1999,harnack2017topological,SMYHDFM:2012,DFHKMMPYS:2013,WPCFMCHMPW:2014,MAC:2014,MW:2014,MAC:2014,YDGS:2015,clark2015spatial,JHHLL:2016,ma2017detection,ma2017detection,Amigo2018,wang2018detecting,pcm}.

The cross map is originated from nonlinear time series analysis~\cite{Takens:1981,PCFS:1980,sauer1991embedology,stark1999delay,stark2003delay,
muldoon1998delay}.
A brief understanding of such a map is as follows. Consider two subsystems: $X$ and $Y$.
In the reconstructed phase space of $X$, if for any state vector at a time a set of neighboring vectors can be found, the set of the cross mapped vectors, which are the partners with equal time of $X$, could be available in $Y$.
The cross map underlying the reconstructed spaces can be written as $Y_t=\Phi(X_t)$ (where $X_t$ and $Y_t$ are delay coordinates with sufficiently large dimensions) for the case of $Y$ {\it unidirectionally} causing $X$ while, mathematically,  its inverse map does not exist~\cite{Amigo2018}.
In practice, using the prior knowledge on the true causality in toy models or/and the assumption on the expanding property of $\Phi$ (representing by its Jacobian's singular value larger than one in the topological causality framework~\cite{harnack2017topological}), scientists developed many practically useful techniques based on the cross map for causality detection. For instance, the ``activity'' method, originally designed to measure the continuity of the inverse of the cross map, compares the divergence of the cross mapped vectors to the state vector in $X$ with the divergence of the independently-selected neighboring vectors to the same state vector{\ }\cite{PRE_Quiroga2000,PhyD_Arnhold1999}. The topological causality measures the divergence rate of the cross mapped vectors from the state vectors in $Y${\ }\cite{harnack2017topological}, and the convergent cross mapping (CCM), increasing the length of time series, compares the true state vector $Y$ with the average of the cross mapped vectors, as the estimation of $Y${\ }\cite{SMYHDFM:2012,DFHKMMPYS:2013,WPCFMCHMPW:2014,MAC:2014,MW:2014,YDGS:2015,clark2015spatial,JHHLL:2016,ma2017detection,Amigo2018,wang2018detecting,PRE_Hirata2010,pcm}.    Then, the change of the divergence or the accuracy of the estimation is statistically evaluated for determining the causation from $Y$ to $X$.   Inversely, the causation from $X$ to $Y$ can be evaluated in an analogous manner.
The above evaluations{\ }\cite{DFHKMMPYS:2013,
WPCFMCHMPW:2014,MAC:2014,MW:2014,YDGS:2015,clark2015spatial,
JHHLL:2016,ma2017detection,harnack2017topological,Amigo2018,
wang2018detecting,PRE_Hirata2010,pcm}
can be understood at a conceptional and qualitative level and perform well in many demonstrations.

In this work, striving for a comprehensive understanding of causal mechanisms and inspired by the cross-map-based techniques,
we develop a mathematically rigorous framework for detecting
causality in nonlinear dynamical systems,
turning eyes towards investigating the original systems from their cross maps,  which is also logically consistent with the natural interpretation of causality as functional dependences{\ }\cite{pearl2009causality,JPearl2009}.
The skills used in cross-map-based methods are assmilated in our framework,
while we directly study the original dynamical systems or the reconstructed systems instead of the cross maps.
The foundation of our framework is the {\it scaling law} for the changing relation of $\varepsilon$ with $\delta$ arising from the continuity for the investigated system, henceforth the term ``continuity scaling''.
In addition to providing a theory, we demonstrate,
using synthetic and real-world data, that our continuity scaling framework is accurate, computationally efficient, widely applicable,
showing advantages over the existing methods.

\section{Continuity Scaling Framework}

\begin{figure*}[t]
\centering
\includegraphics[width=0.7\textwidth]{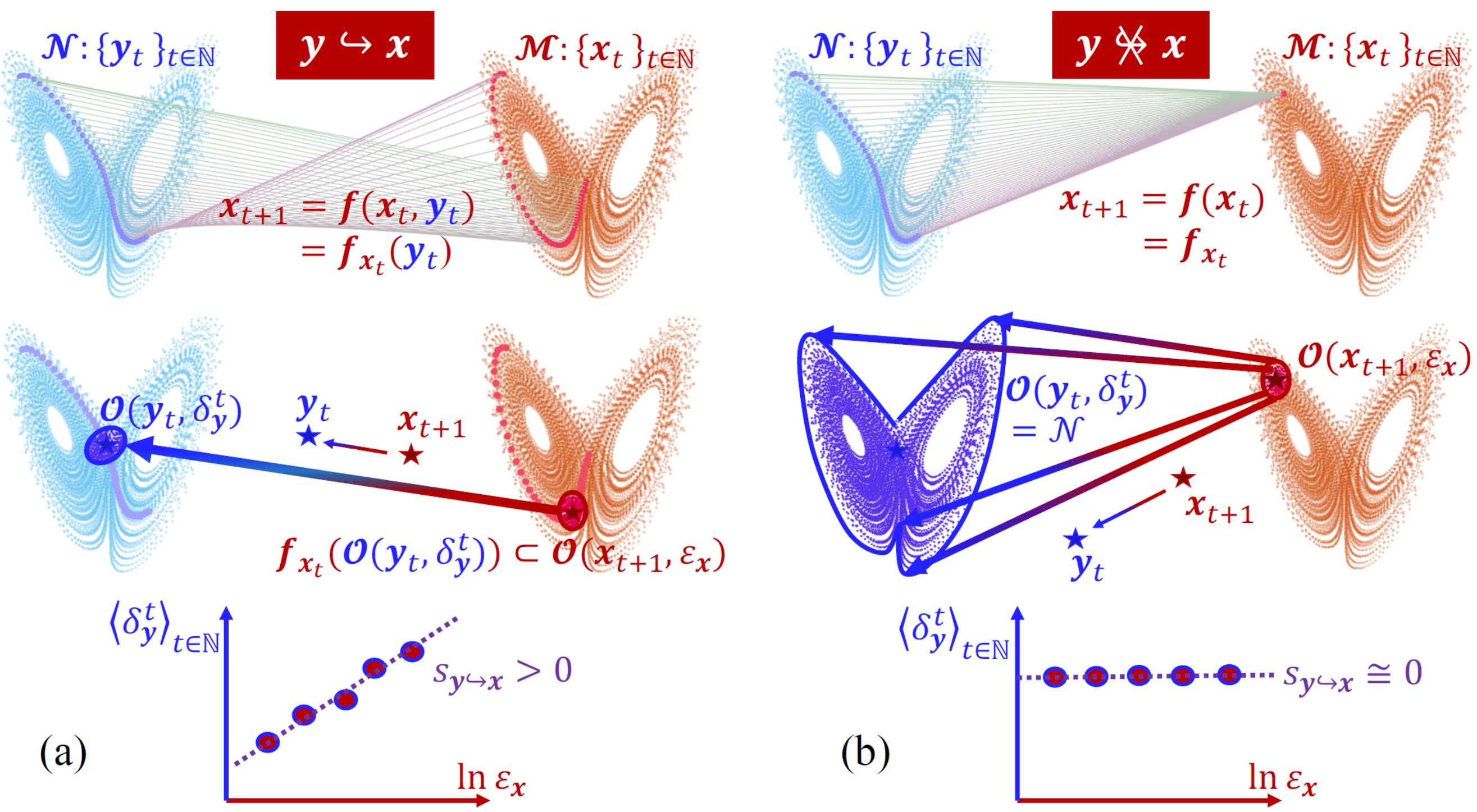}
\subfigure{\label{fig1a}}
\subfigure{\label{fig1b}}
\caption{ Illustration of causal relation between two sets of
dynamical variables. (a) Existence of causation from $\bm{y}$ in $\mathcal{N}$ to
$\bm{x}$ in $\mathcal{M}$, where each correspondence from $\bm{x}_{t+1}$ to $\bm{y}_t$ is one-to-one, represented by the line or the arrow, respectively, in the upper and the middle panels.
In this case, a change in $\ln{\varepsilon_{\bm{x}}}$ results in
a direct change in $\delta_{\bm{y}}$ (the lower panel) with $\varepsilon_{\bm{x}}$ and
$\delta_{\bm{y}}$ denoting the neighborhood size of the resulting
variable $\bm{x}$ and of the causal variable $\bm{y}$, respectively. (b) Absence of causation from $\bm{y}$ to $\bm{x}$, where every point on each trajectory, $\{\bm{y}_{t}\}$, in
$\mathcal{N}$ could be the correspondent point from $\bm{x}_{t+1}$ in $\mathcal{M}$ (the upper panel) and thus every point in
$\mathcal{N}$ belongs to the largest $\delta$-neighborhood of
$\bm{y}_{t}$ (the middle panel). In this case, $\delta_{\bm{y}}$ does not depend on $\varepsilon_{\bm{x}}$ (the lower panel). Also refer to the supplemental animation for illustration.
}\label{fig1}
\end{figure*}

To explain the mathematical idea behind the development of our framework, we
use the following class of discrete time dynamical systems:
$\bm x_{t+1}=\bm f(\bm x_t,\bm y_t)$ and
$\bm y_{t+1}=\bm g(\bm x_t,\bm y_t)$ for $t\in\mathbb{N}$,
where the state variables $\{\bm x_t\}_{t\in\mathbb{N}}$, $\{\bm y_t\}_{t\in\mathbb{N}}$ evolve in the compact manifolds $\mathcal{M}$, $\mathcal{N}$ of dimension $D_{\mathcal{M}}$, $D_{\mathcal{N}}$ under
sufficiently smooth map $\bm{f}$, $\bm{g}$, respectively. We adopt the common recognition of causality in dynamical systems.
\begin{shaded}
\noindent{\bf Definition.} If the dependence of
$\bm{f}(\bm{x},\bm{y})$ on $\bm{y}$ is nontrivial (i.e., a directional coupling exists), a variation
in $\bm{y}$ will result in a change in the value of
$\bm{f}(\bm{x},\bm{y})$ for any given $\bm{x}$, which, according to
the natural interpretation of
causality~\cite{spirtes2000causation,pearl2009causality}, admits that $\bm{y}:\{\bm y_t\}_{t\in\mathbb{N}}$ has a direct causal effect on
$\bm{x}:\{\bm x_t\}_{t\in\mathbb{N}}$, denoted by
$\bm{y}\hookrightarrow\bm{x}$, as shown in the upper panel of Fig.~\ref{fig1a}.
\end{shaded}

We now interpret the causal relationship stipulated by the
continuity of a function. Let
$\bm f_{\bm x_{\text{g}}}(\cdot)\triangleq\bm f(\bm x_{\text{g}},\cdot)$ for
a given point $\bm x_{\text{g}}\in\mathcal{M}$. For any
$\bm y_P\in \mathcal{N}$, we denote its image under the given function by
$\bm x_I\triangleq\bm f_{\bm x_{\text{g}}}(\bm y_P)$.
Applying the logic statement of a continuous function to
$\bm f_{\bm x_{\text{g}}}(\cdot)$,  we have that,
\textit{for any neighborhood $\mathcal{O}(\bm x_I,\varepsilon)$ centered at
$\bm x_I$ and of radius $\varepsilon>0$, there exists a neighborhood
$\mathcal{O}(\bm y_P,\delta)$ centered at $\bm y_P$ of radius $\delta>0$,
such that
$\bm f_{\bm x_{\text{g}}}(\mathcal{O}(\bm y_P,\delta))\subset \mathcal{O}(\bm x_I,\varepsilon)$.} The neighborhood and its radius are defined by
$$
   \mathcal{O}(\bm{p},h)=\big\{ \bm{s}\in\mathcal{S}~\big|~
   \text{dist}_{\mathcal{S}}(\bm{s},\bm{p})<h
   \big\},~\bm{p}\in\mathcal{S},~h>0,
$$
where $\text{dist}_{\mathcal{S}}(\cdot,\cdot)$ represents an appropriate
metric describing the distance between two given points in a specified manifold
$\mathcal{S}$ with $\mathcal{S}=\mathcal{M}$ or $\mathcal{N}$. The meaning of
this mathematical statement is that, if we have a neighborhood of the resulting
variable $\bm x_{I}$ first, we can then find a neighborhood for the causal variable
$\bm y_{P}$ to satisfy the above mapping and inclusion relation. This operation
of ``first-$\varepsilon$-then-$\delta$'' provides a rigorous base for the
principle that the information about the resulting variable can be used to estimate the information of the causal variable and therefore to ascertain
causation, as indicated by the long arrow in the middle panels of
Fig.~\ref{fig1a}. Note that, the existence of the $\delta>0$ neighborhood is
always guaranteed for a continuous map $\bm f_{\bm x_{\text{g}}}$.
In fact, due to the compactness of the manifold $\mathcal{N}$, a largest
value of $\delta$ exists. However, if $\bm{y}_{P}$ does not have an explicit
causal effect on the variable $\bm{x}_{I}$, i.e., $\bm{f}_{\bm{x}_{\rm g}}$
is independent of $\bm{y}_{P}$, the existence of $\delta$ is still assured
but it is independent of the value of $\varepsilon$, as shown in the upper
panel of Fig.~\ref{fig1b}. This means that merely determining the existence
of a $\delta$-neighborhood is not enough for inferring causation - it is
necessary to vary $\varepsilon$ systematically and to examine the scaling
relation between $\delta$ and $\varepsilon$. In the following we discuss
a number of scenarios.
\\\\
{\bf Case I}: Dynamical variables $\{(\bm x_t,\bm y_t)\}_{t\in\mathbb{N}}$
are fully measurable. For any given constant $\varepsilon_{\bm x}>0$, the set
$\{\bm x_{\tau}\in\mathcal{M}~|~\tau\in I^t_{\bm x}(\varepsilon_{\bm x})\}$
can be used to approximate the neighborhood
$\mathcal{O}(\bm x_{t+1},\varepsilon_{\bm{x}})$, where the time index set is
\begin{equation}\label{indexset}
    I^t_{\bm x}(\varepsilon_{\bm x})
\triangleq \left\{\tau\in\mathbb{N}~|~\text{dist}_{\mathcal{M}}(\bm x_{t+1},\bm x_{\tau})<\varepsilon_{\bm x}\right\}.
\end{equation}
The radius $\delta^t_{\bm y}=\delta^t_{\bm y}(\varepsilon_{\bm x})$
of the neighborhood $\mathcal{O}(\bm y_t,\delta^t_{\bm y})$ satisfying
$\bm f_{\bm{x}_{\rm g}=\bm{x}_t}(\mathcal{O}(\bm y_t,\delta^t_{\bm y}))\subset \mathcal{O}(\bm x_{t+1},\varepsilon_{\bm{x}})$
can be estimated as
\begin{equation}\label{delta}
\delta^t_{\bm y}(\varepsilon_{\bm x})
  \triangleq \left\{\#[ \bar{I}^t_{\bm x}(\varepsilon_{\bm x})]\right\}^{-1}\sum_{\tau\in \bar{I}^t_{\bm x}(\varepsilon_{\bm x})} \text{dist}_{\mathcal{N}}(\bm y_t,\bm y_{\tau-1}),
\end{equation}
where $\#[\cdot]$ is the cardinality of the given set and the index set is
$\bar{I}^t_{\bm x}(\varepsilon_{\bm x})\triangleq \{\tau\in I^t_{\bm x}(\varepsilon_{\bm x})~|~ ~\text{dist}_{\mathcal{M}}(\bm x_{t},\bm x_{\tau-1})<\varepsilon_{\bm x}\}$.

The strict mathematical steps for estimating $\delta^t_{\bm y}$ are given in Section II of Supplementary Information (SI). We emphasize that here correspondence between $\bm{x}_{t+1}$ and $\bm{y}_t$ is investigated, differing from the cross-map-based methods, with one-step time difference naturally arising. This consideration yields a key condition [DD], which is only need when considering the original iteration/flow and whose detailed description and universality are demonstrated in SI. We reveal a linear scaling law between $\langle\delta^t_{\bm y}\rangle_{t\in\mathbb{N}}$ and
$\ln{\varepsilon_{\bm x}}$, as shown in the lower panels of Fig.~\ref{fig1},
whose slope $s_{\bm{y}\hookrightarrow\bm{x}}$ is an indicator of
the correspondent relation between $\varepsilon$ and $\delta$ and hence
the causal relation $\bm{y}\hookrightarrow\bm{x}$. Here, $\langle\cdot\rangle_{t\in\mathbb{N}}$ denotes the average over time. In particular, a larger
slope value implies a stronger causation in the direction from $\bm y$ to
$\bm x$ as represented by the map functions $\bm{f}(\bm{x}_t,\bm{y}_t)$
[Fig.~\ref{fig1a}], while a near zero slope indicates null causation
in this direction [Fig.~\ref{fig1b}]. Likewise, possible causation in the
reversed direction, $\bm{x}\hookrightarrow\bm{y}$, as represented by the
function $\bm{g}(\bm{x}_t,\bm{y}_t)$, can be assessed analogously.
And the unidirectional case when $\bm{f}(\bm{x},\bm{y})=\bm{f}_{0}(\bm{x})$ independent of $\bm{y}$ is uniformly considered in Case II.
We summarize the consideration below and an argument for the generic existence of the scaling law is provided in Section II of SI.
\begin{shaded}
\noindent{\bf Theorem.} For dynamical variables $\{(\bm x_t,\bm y_t)\}_{t\in\mathbb{N}}$ measured directly from the dynamical systems, if the slope $s_{\bm{y}\hookrightarrow\bm{x}}$ defined above is zero, no causation exists from $\bm{y}$ to $\bm{x}$. Otherwise, a directional coupling can be confirmed from $\bm{y}$ to $\bm{x}$ and the slope $s_{\bm{y}\hookrightarrow\bm{x}}$ increases monotonically with the coupling strength.
\end{shaded}

\begin{figure*}[t]
\centering
\includegraphics[width=0.6\textwidth]{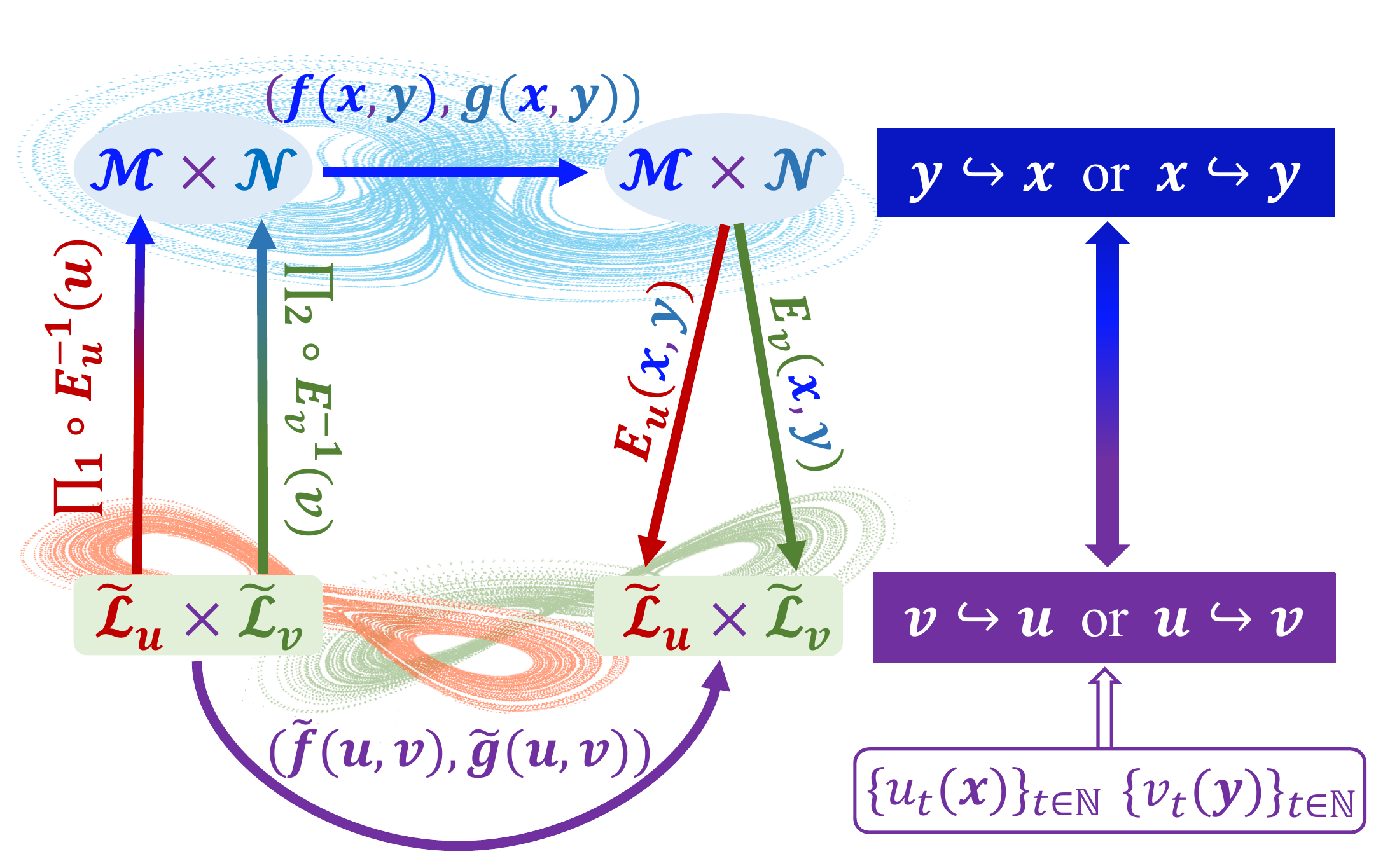}
\caption{Illustration of system dynamics before and after embedding for Case II. In the left panel, the arrows describe how the original systems $(\bm{f},\bm{g})$ is equivalent to the system $(\tilde{\bm{f}},\tilde{\bm{g}})$ after embedding.  In the right panel, causation
between the internal variables $\bm{x}$ and $\bm{y}$ can be ascertained by
detecting the causation between the variables $\bm{u}$ and $\bm{v}$
reconstructed from measured time series.}
\label{fig2}
\end{figure*}

\noindent{\bf Case II}: The dynamical variables $\{(\bm x_t,\bm y_t)\}_{t\in\mathbb{N}}$
are not directly accessible but measurable time series
$\{u_t\}_{t\in\mathbb{N}}$ and $\{v_t\}_{t\in\mathbb{N}}$ are available, where
$u_t=u(\bm x_t)$ and $v_t=v(\bm y_t)$ with $u$:
$\mathcal{M}\to\mathbb{R}^{r_u}$ and $v$: $\mathcal{N}\to\mathbb{R}^{r_v}$
being smooth observational functions. To assess causation from $\bm y$ to
$\bm x$, we assume one-dimensional observational time series (for simplicity):
$r_u=r_v=1$, and use the classical delay-coordinate embedding
method~\cite{Takens:1981,PCFS:1980,sauer1991embedology,stark1999delay,
stark2003delay,Cummins2015,muldoon1998delay} to reconstruct the phase space:
$\bm u_t=(u_t,u_{t+\tau_u},\cdots,u_{t+(d_u-1)\tau_u})^{\top}$ and
$\bm v_t=(v_t,v_{t+\tau_v},\cdots,v_{t+(d_v-1)\tau_v})^{\top}$,
where $\tau_{u,v}$ is the delay time and
$d_{u,v}>2(D_{\mathcal{M}}+D_{\mathcal{N}})$ is the embedding dimension
that can be determined using some standard criteria~\cite{kantz2004nonlinear}.
As illustrated in Fig.~\ref{fig2}, the dynamical evolution of the reconstructed
states $\{(\bm u_t,\bm v_t)\}_{t\in\mathbb{N}}$ is governed by
\begin{equation}\label{equ:obsys}
  \bm u_{t+1}=\tilde{\bm f}(\bm u_t,\bm v_t), ~~
  \bm v_{t+1}=\tilde{\bm g}(\bm u_t,\bm v_t).
\end{equation}
The map functions can be calculated as
$\tilde{\bm f}(\bm{u},\bm{v}) \triangleq
\bm E_u \circ [\bm f,\bm g]\left(\bm \Pi_1\circ\bm E_u^{-1}(\bm{u}),\bm \Pi_2\circ\bm E_v^{-1}(\bm{v})\right)$,
$\tilde{\bm g}(\bm{u},\bm{v}) \triangleq
\bm E_v \circ [\bm f,\bm g]\left(\bm \Pi_1\circ\bm E_u^{-1}(\bm{u}),\bm \Pi_2\circ\bm E_v^{-1}(\bm{v})\right)$,
where the embedding (diffeomorphism) $\bm E_{s}$: $\mathcal{M}\times\mathcal{N}\to \widetilde{\mathcal{L}}_s\subset\mathbb{R}^{d_s}$ with $\widetilde{\mathcal{L}}_s\triangleq\bm E_{s}(\mathcal{M}\times\mathcal{N})$, $s = u \ \mathrm{or} \ v$, is given by
$$
\begin{array}{l}
\bm E_{u}(\bm x,\bm y)\triangleq
\big(u(\bm x),u\circ\bm{\Pi}_1\circ[\bm f,\bm g]^{\tau_u}(\bm x,\bm y), u\circ\bm{\Pi}_1\circ\\
\hspace{1cm} [\bm f,\bm g]^{2\tau_u}(\bm x,\bm y),  \cdots, u\circ\bm{\Pi}_1\circ[\bm f,\bm g]^{(d_u-1)\tau_u}(\bm x,\bm y)
\big), \\
\bm E_{v}(\bm x,\bm y)\triangleq
\big(v(\bm y),v\circ\bm{\Pi}_2\circ [\bm f,\bm g]^{\tau_v}(\bm x,\bm y), v\circ\bm{\Pi}_2\circ\\
\hspace{1cm}[\bm f,\bm g]^{2\tau_v}(\bm x,\bm y),
\cdots,  v\circ\bm{\Pi}_2\circ[\bm f,\bm g]^{(d_v-1)\tau_v}(\bm x,\bm y)
\big),
\end{array}
$$
with the inverse function $\bm E_{s}^{-1}$ defined on
$\widetilde{\mathcal{L}}_s$, $[\bm{f},\bm{g}]^{k}$ representing the $k$th
iteration of the map, and the projection mappings as
$\bm \Pi_1(\bm x,\bm y)=\bm x$ and $\bm \Pi_2(\bm x,\bm y)=\bm y$. Case II
has now been reduced to Case I and our continuity scaling framework can
be used to ascertain the causation from $\bm{v}$ to $\bm{u}$ based on the
measured time series with the indices $I^t_{\bm{u}}(\varepsilon_{\bm{u}})$,
$\delta^t_{\bm{v}}(\varepsilon_{\bm{u}})$ and $s_{\bm{v}\hookrightarrow\bm{u}}$ [Eqs.~\eqref{indexset} and \eqref{delta}].

Does the causation from $\bm{v}$ to $\bm{u}$ imply causation from $\bm{y}$
to $\bm{x}$? The answer is affirmative, which can be argued, as follows.
If the original map function $\bm{f}$ is independent of $\bm{y}$:
$\bm f(\bm x,\bm y)=\bm f_0(\bm x)$, there is no causation from $\bm{y}$ to
$\bm{x}$. In this case, the embedding $\bm E_u(\bm x,\bm y)$ becomes
independent of $\bm{y}$, degenerating into the form of
$\bm E_u(\bm x,\bm y)=\bm E_{u0}(\bm x)$, a diffeomorphism from $\mathcal{M}$
to $\widetilde{\mathcal{L}}_{u0}=\bm E_{u0}(\mathcal{M})$ only. As a result,
Eq.~\eqref{equ:obsys} becomes: $\bm{u}_{t+1}=\tilde{\bm f}_0(\bm{u}_t)$ and
$\bm{v}_{t+1}=\tilde{\bm g}(\bm{u}_t,\bm{v}_t)$, where
$\tilde{\bm f}_0(\bm{u})=\bm{E}_{u0}\circ\bm{f}\circ\bm{E}_{u0}^{-1}(\bm{u})$
and the resulting mapping $\tilde{\bm{f}}_{0}$ is independent of $\bm{v}$.
The independence can be validated by computing the slope
$s_{\bm{v}\hookrightarrow\bm{u}}$ associated with the scaling relation
between $\langle\delta^t_{\bm v}\rangle_{t\in\mathbb{N}}$ and
$\ln{\varepsilon_{\bm{u}}}$, where a zero slope indicates null causation
from $\bm{v}$ to $\bm{u}$ and hence null causation from $\bm{y}$
to $\bm{x}$. Conversely, a finite slope signifies causation between the
variables. Thus, any type of causal relation (unidirectional or
bi-directional) detected between the reconstructed state variables
$\{(\bm{u}_{t},\bm{v}_{t})\}_{t\in\mathbb{N}}$ implies the same type of causal
relation between the internal but inaccessible variables $\bm{x}$ and $\bm{y}$
of the original system.
\\\\
{\bf Case III}: The structure of the internal variables is completely unknown.
Given the observational functions
$\tilde{u},\tilde{v}$: $\mathcal{M}$$\times$$\mathcal{N}$$\to$$\mathbb{R}$ with
$\tilde{u}_t=\tilde{u}(\bm{x}_t,\bm{y}_{t})$ and
$\tilde{v}_t=\tilde{v}(\bm{x}_t,\bm{y}_{t})$, we first reconstruct the state
space:
$\tilde{\bm u}_t=(\tilde{u}_t, \tilde{u}_{t+\tau},\cdots,\tilde{u}_{t+(d-1)\tau})^{\top}$ and
$\tilde{\bm v}_t=(\tilde{v}_t, \tilde{v}_{t+\tau},\cdots,\tilde{v}_{t+(d-1)\tau})^{\top}$.
To detect and quantify causation from $\tilde{\bm{v}}$ to $\tilde{\bm{u}}$
(or vice versa), we carry out a continuity scaling analysis with the modified
indices
$I^t_{\tilde{\bm{u}}}(\varepsilon_{\tilde{\bm{u}}})$, $\delta^t_{\tilde{\bm{v}}}(\varepsilon_{\tilde{\bm{u}}})$ and $s_{\tilde{\bm{v}}\hookrightarrow\tilde{\bm{u}}}$.
Differing from Case II, here, due to the lack of knowledge
about the correspondence structure between the internal and observational
variables, a causal relation for the latter does not definitely imply the same
for the former.
\\\\
{\bf Case IV}: Continuous-time dynamical systems possessing a sufficiently
smooth flow $\{\bm S_t;t\in\mathbb{R}\}$ on a compact manifold $\mathcal{H}$:
$\text{d}\bm S _t(\bm u_0)/\text{d}t = \bm \chi(\bm S_t(\bm u_0))$, where
$\bm{\chi}$ is the vector field. Let
$\{\hat{u}_{t=\omega n+\nu}\}_{n\in\mathbb{Z}}$ and
$\{\hat{v}_{t=\omega n+\nu}\}_{n\in\mathbb{Z}}$ be two respective time series
from the smooth observational functions
$\hat{u},~\hat{v}$: $\mathcal{H}$$\rightarrow$$\mathbb{R}$ with
$\hat{u}_t=\hat{u}(\bm{S}_t)$ and $\hat{v}_t=\hat{v}(\bm{S}_t)$,  where $
1/\omega$ is the sampling rate and $\nu$ is the time shift. Defining
$\bm{\Xi}\triangleq\bm{S}_{\omega}$: $\mathcal{H}$$\rightarrow$$\mathcal{H}$
and $\hat{\bm{S}}_n\triangleq {\bm{S}_{\omega n+\nu}}(\bm{u}_{0})$, we
obtain a discrete-time system as $\hat{\bm S}_{n+1}=\bm \Xi (\hat{\bm S}_n)$
with the observational functions as
$\hat{u}_n=\hat{u}(\hat{\bm{S}}_n)$ and $\hat{v}_n=\hat{v}(\hat{\bm{S}}_n)$,
reducing the case to Case III and rendering applicable our continuity scaling
analysis to unveil and quantify the causal relation between
$\{\hat{u}_{t=\omega n+\nu}\}_{n\in\mathbb{Z}}$ and
$\{\hat{v}_{t=\omega n+\nu}\}_{n\in\mathbb{Z}}$. If the domains of
$\hat{u}$ and $\hat{v}$ have their own restrictions on some particular
subspaces, e.g., $\hat{u}$: $\mathcal{H}_{u}$$\rightarrow$$\mathbb{R}$ and
$\hat{v}$: $\mathcal{H}_{v}$$\rightarrow$$\mathbb{R}$ with
$\mathcal{H}=\mathcal{H}_{u}\oplus\mathcal{H}_{v}$, the case is further
reduced to Case II, so the detected causal relation between the observational
variables imply causation between the internal variables belonging to their
respective subspaces.

\begin{figure*} [ht!]
\centering
{\includegraphics[width=0.7\textwidth]{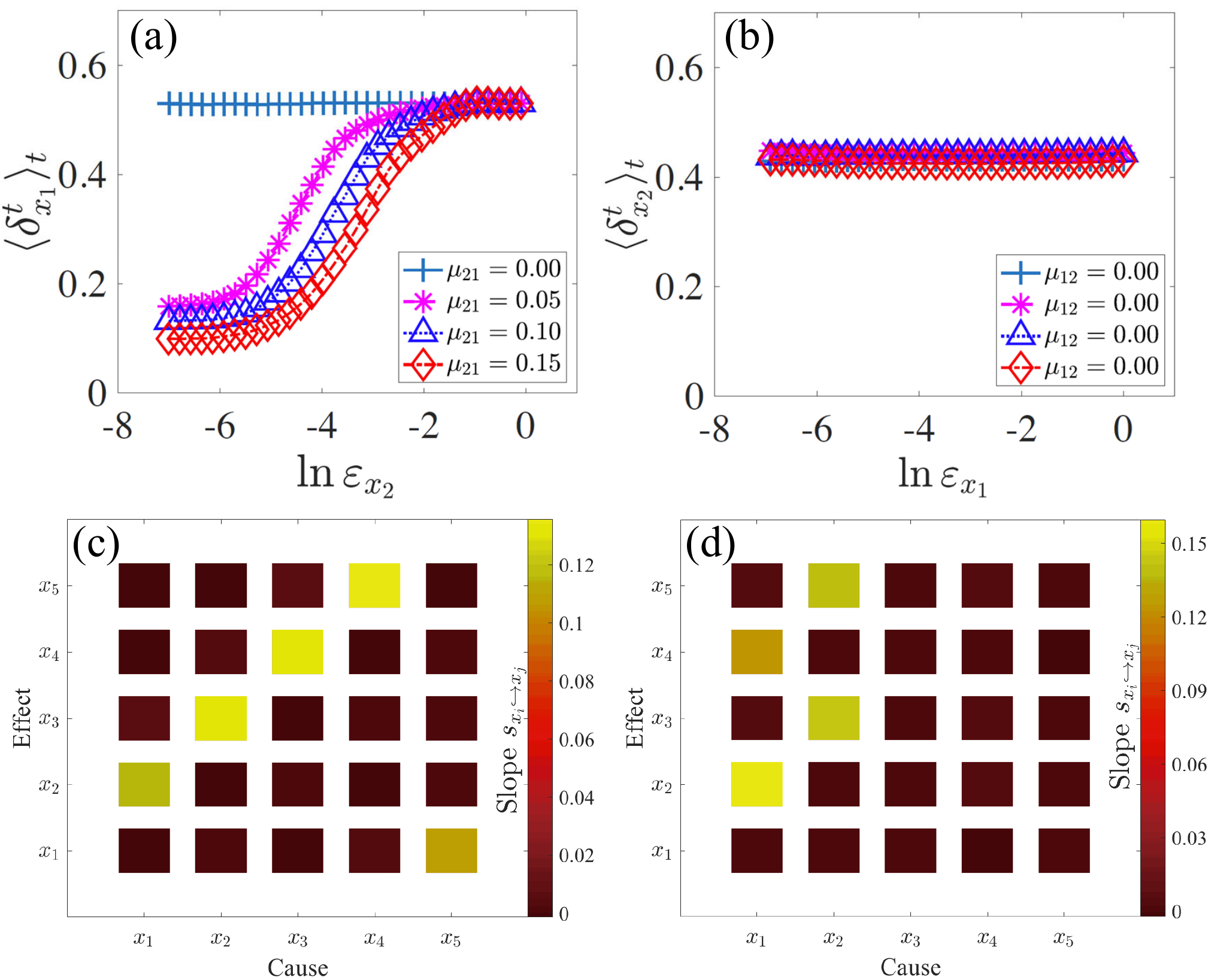}}
\subfigure{\label{fig3a}}
\subfigure{\label{fig3b}}
\subfigure{\label{fig3c}}
\subfigure{\label{fig3d}}
\caption{\label{fig3}
Ascertaining and characterizing causation in various ecological
systems of interacting species. (a,b) Unidirectional causation of two coupled
species. In (a), the values of the slope $s_{x_{1}\hookrightarrow x_{2}}$
associated with the causal relation $x_{1}\hookrightarrow x_{2}$ are
approximately 0.0004, 0.1167, 0.1203, and 0.1238 for four different values
of the coupling parameter $\mu_{21}$. (b) Near zero slope values
$s_{x_{2}\hookrightarrow x_{1}}$ for $x_{2}\hookrightarrow x_{1}$, indicating
its nonexistence. (c,d) Inferred causal network of five species whose
interacting structure is, respectively, that of a ring:
$x_{i}\hookrightarrow x_{i+1(\mathrm{mod}~5)}$ ($i=1,\cdots,5$) and
of a tree: $x_{j}\hookrightarrow x_{j+1,j+3}$ ($j=1,2$), where the
estimated slope values are color-coded.
Results of a statistical analysis of the accuracy and reliability of the determined causal interactions are presented in SI section III.
Time series of length $5000$ are used in all these simulations.
The embedding parameters are $\tau_{s}=1$ and
$d_{s}=3$ with $s=x_{1},\cdots,x_{5}$.
}
\end{figure*}

\section{Demonstrations: From Complex Dynamical Models to Real-World Networks}

To demonstrate the efficacy of our continuity scaling framework and its
superior performance, we have carried out extensive numerical tests with a
large number of synthetic and empirical datasets, including those from gene
regulatory networks as well as those of air pollution and hospital admission.
The practical steps of the continuity scaling framework together with the significance test procedures are described in Methods. We present three representative examples here, while leaving others of significance to SI.

The first example is an ecological model of two unidirectionally interacting
species: $x_{1,t+1}=x_{1,t}(3.8-3.8 x_{1,t}-\mu_{12} x_{2,t})$ and
$x_{2,t+1}=x_{2,t}(3.7-3.7 x_{2,t}-\mu_{21} x_{1,t})$. With time series
$\{(x_{1,t},x_{2,t})\}_{t\in\mathbb{N}}$ obtained from different values
of the coupling parameters, our continuity scaling framework yields correct results
of different degree of unidirectional causation, as shown in
Figs.~\ref{fig3a}-\ref{fig3b}. In all cases, there exists a reasonable
range of $\ln{\varepsilon_{x_2}}$ (neither too small nor too large) from which
the slope $s_{x_1\hookrightarrow x_2}$ of the linear scaling can be extracted.
The statistical significance of the estimated slope values and consequently
the strength of causation can be assessed with the standard $p$-value
test~\cite{lancaster2018surrogate} (Methods and SI). An ecological model with
bidirectional coupling has also been tested (see Section III of SI). Figures~\ref{fig3c}-\ref{fig3d} show the results from ecological networks
of five mutually interacting species on a ring and on a tree structure, respectively, where the color-coded slope values reflect accurately the interaction
patterns in both cases.

\begin{figure*}[t]
\centering
{\includegraphics[width=0.7\textwidth]{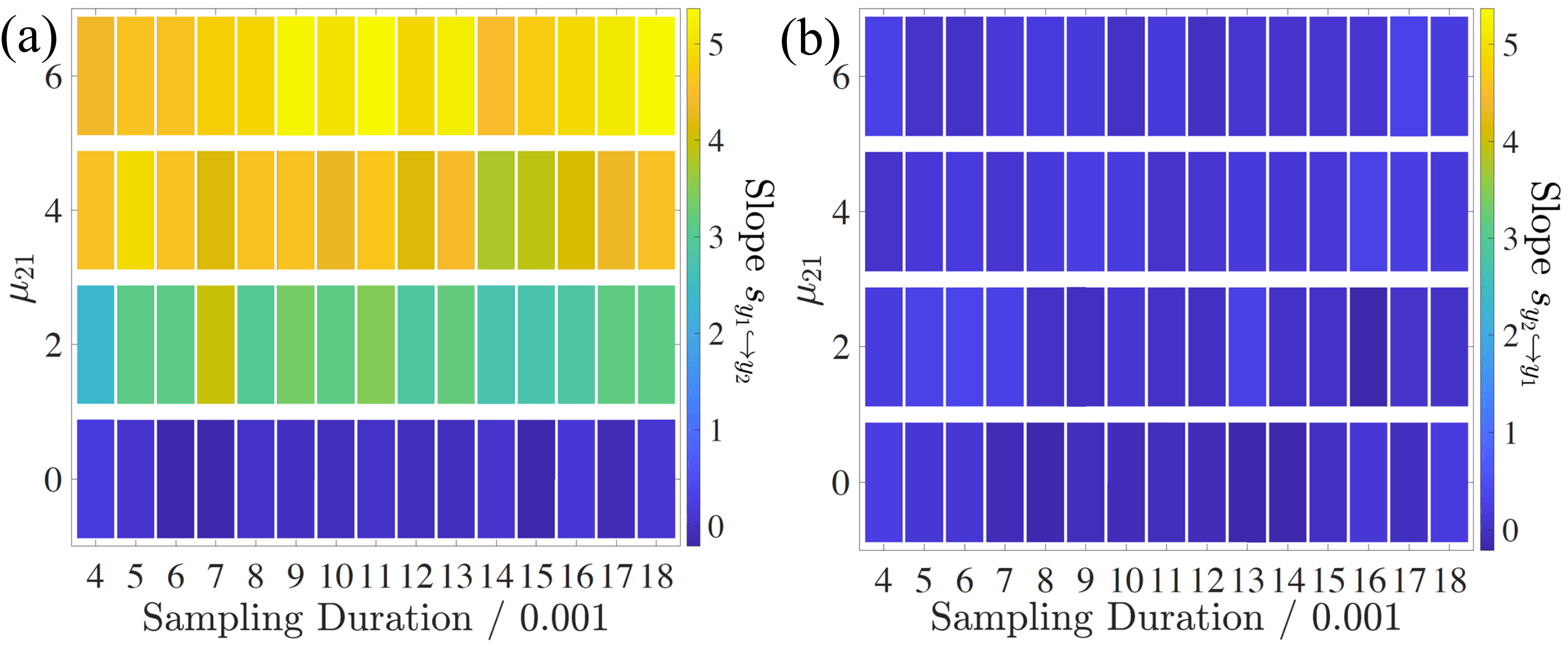}}
\subfigure{\label{fig4a}}
\subfigure{\label{fig4b}}
\caption{\label{fig4}
Detecting causation in the unidirectionally coupled Lorenz system.
The results are for different values of $\mu_{21}$ ($\mu_{12}=0$) and sampling
rate $1/\omega$. (a,b) Color-coded values of the slopes
$s_{y_{1}\hookrightarrow y_{2}}$ and $s_{y_{2}\hookrightarrow y_{1}}$,
respectively. The integration time step is $10^{-3}$ and the embedding parameters are $d_{s}=7$, $\tau_{s}\approx0.05$ with
$\omega|\tau_{s}$ ($s=y_{1}$ or $y_{2}$).
See Section III and Tab. S9 of SI for all the other parameters including the time series lengths used in the simulations.
}
\end{figure*}

The second example is the coupled Lorenz system:
$\dot x_i=\sigma_i(y_i-x_i)+\mu_{ij}x_j$, $\dot y_i=x_i(\rho_i-z_i)-y_i$,
$\dot z_i=x_i y_i -\beta_i z_i$ with $i,j=1,2$ and $i\not=j$. We use
time series $\{y_{1,t}, y_{2,t}\}_{t=n\omega}$ for detecting different
configurations of causation (see Section III of SI). Figure~\ref{fig4} presents the
overall result, where the color-coded estimated values of the slope from
the continuity scaling are shown for different combinations of the sampling
rate $1/\omega$ and coupling strength. Even with relatively low sampling rate,
our continuity scaling framework can successfully detect and quantify the
strength of causation.
Note that the accuracy does not vary monotonously with the sampling rate,
indicating the potential of our framework to ascertain and quantify causation even with rare data. Moreover, the proposed index can accurately reflect the true causal strength (denoting by the coupling parameter),
which is also evidenced by numerical tests in Sections III and IV of SI.
Robustness tests against different noise perturbations are provided in Section III of SI demonstrating the practical usefulness of our framework.
Additionally, analogous to the first example, we
present in SI several examples on causation detection in the coupled Lorenz system with nonlinear couplings,
and the R\"ossler-Lorenz system, etc.,
which further demonstrates the generic efficacy of our framework.

\begin{figure*}[t]
\centering
{\includegraphics[width=0.7\textwidth]{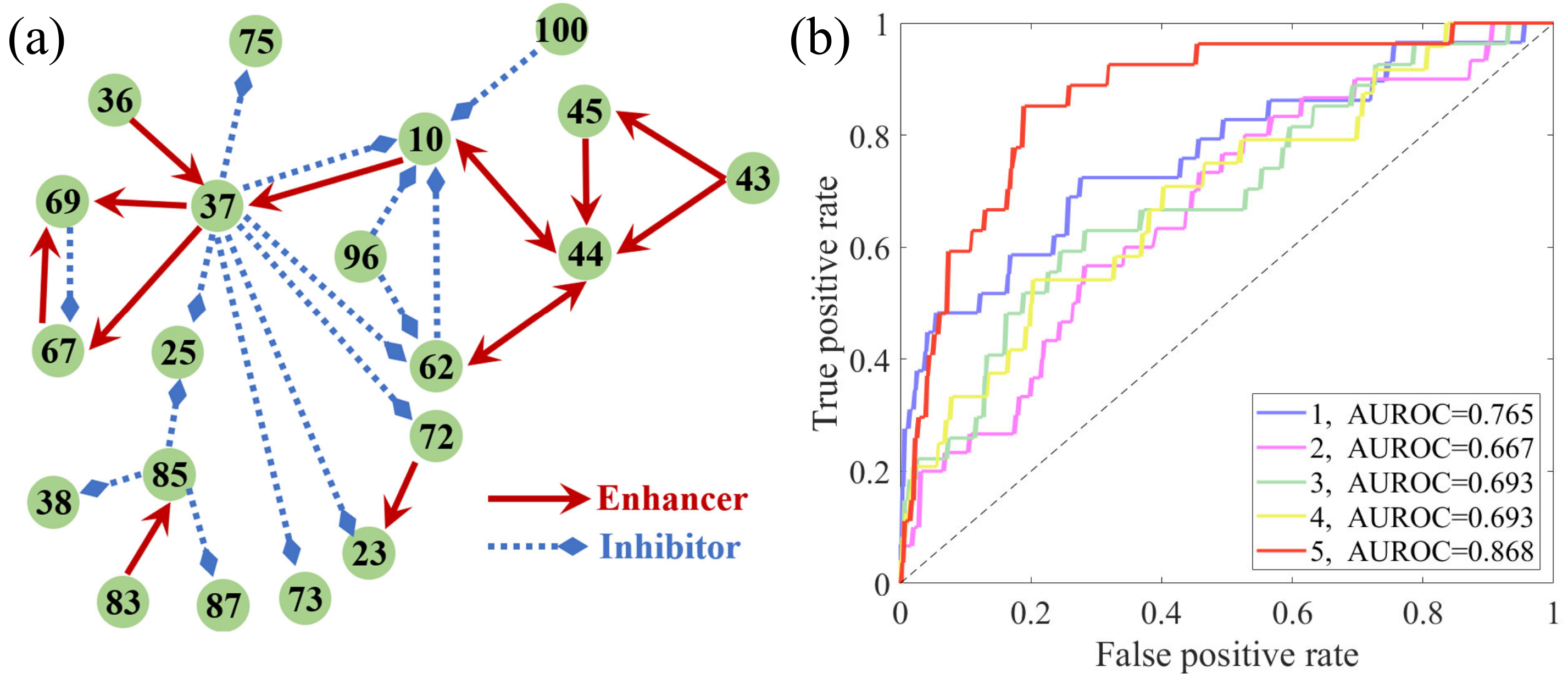}}
\subfigure{\label{fig5a}}
\subfigure{\label{fig5b}}
\caption{\label{fig5}
Detecting causal interactions in five GRNs. (a) One representative GRN containing 20 randomly-selected genes. Other four structures can be find in Fig.~S17 of SI. (b) The ROC curves as well as their AUROC values demonstrate the efficacy of our framework.}
\end{figure*}

In addition, we present study on several real-world dataset,
which brings new insights to the evolutionary mechanism of underlying systems.
We study gene expression data from DREAM4 {\it in silico} Network Challenge~\cite{marbach2009generating,marbach2010revealing}, whose intrinsic gene regulatory networks (GRNs) are known for verification [Fig.~\ref{fig5a} and Fig.~S17 of SI].
Applying our framework to these data, we ascertain the causations between each pair of genes by using the continuity scaling framework.  The corresponding ROC curves for five different networks as well as their AUROC values are shown in Fig.~\ref{fig5b}, which indicates a high detection accuracy in dealing with real-world data.

\begin{figure*}[t]
\centering
\includegraphics[width=0.7\textwidth]{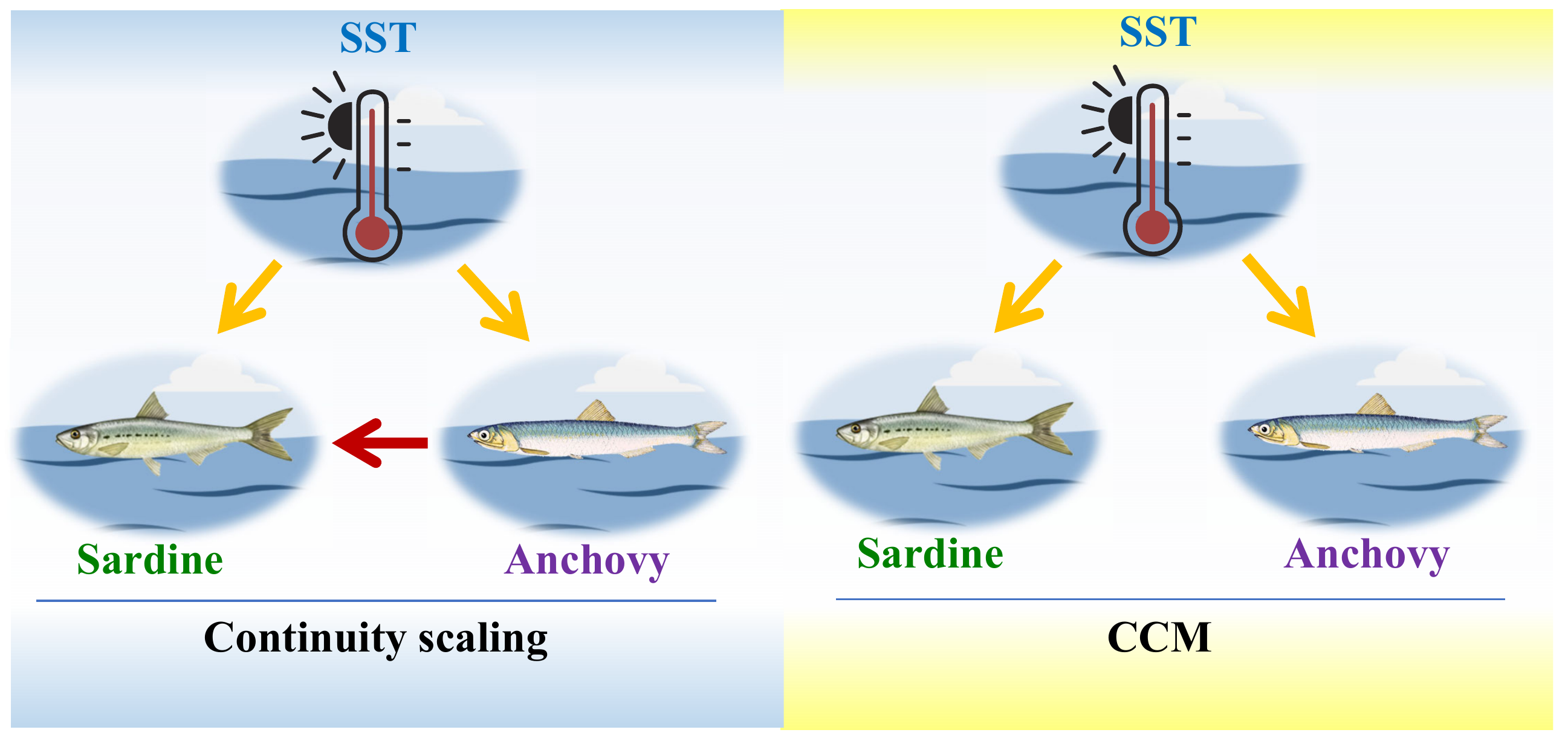}
\caption{\label{E25mainpic}
The comparison of  causal network structure detected by continuity scaling and CCM among sea surface temperature, sardine and anchovy.}
\end{figure*}

We then test the causal relationship in a marine ecosystem consisting of Pacific sardine landings, northern anchovy landings and sea surfae temperature (SST).
We reveal new findings to support the competing relationship hypothesis stated in{\ }\cite{Lasker1983Ocean} which cannot be detected by CCM{\ }\cite{SMYHDFM:2012}.
As pointed out in Fig. \ref{E25mainpic}, while common influence from SST to both species is verified with both methods, our continuity scaling additionally illuminates notable influence from anchovy to sardine
with its reverse direction being less significant.
While competing relationship plays an important role in ecosystems, continuity scaling can reveal more essential interaction mechanism.
See Section III.E  of SI for more details.

\begin{figure*}[ht!]
\centering
\includegraphics[width=0.98\textwidth]{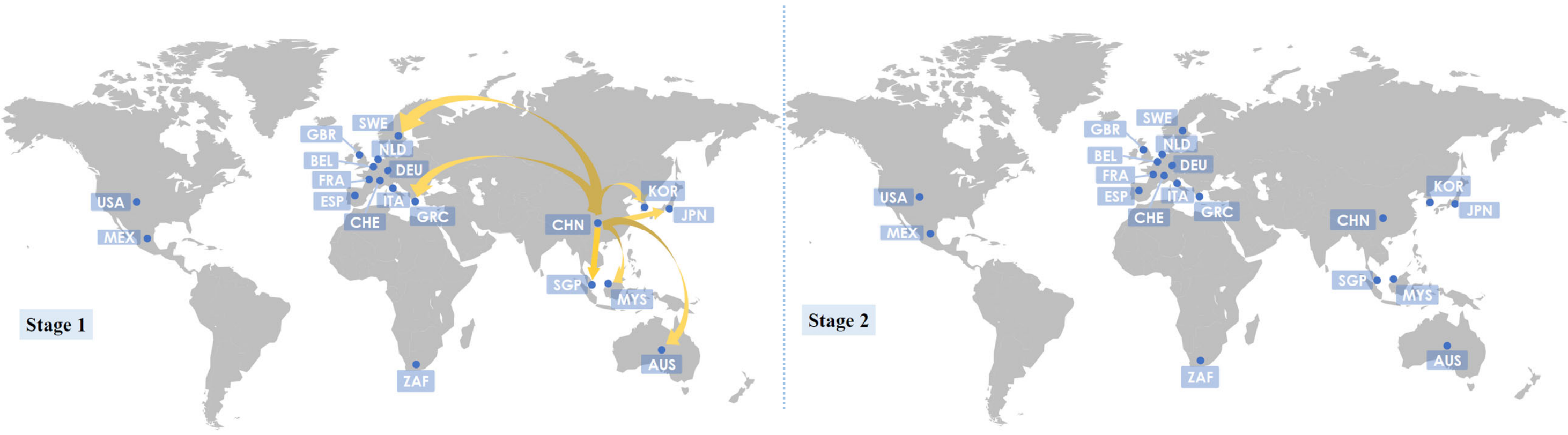}
\caption{\label{E24mainpic}
The  causal effect from China to other countries of the COVID-19 pandemic detected by continuity scaling between Stage $1$ and $2$.
Here, Stage $1$ is from January 22$^{\text{nd}}$ 2020 to April 30$^{\text{th}}$ 2020,
and Stage $2$ is from May 1${}^{\text{st}}$ 2020
to February 15${}^{\text{th}}$ 2021.
For those detected causal links between all countries, refer to Section III.E and Fig.~S18 of SI.
These maps are for illustration only.
}
\end{figure*}

Moreover, we study the transmission mechanism of the recent COVID-19 pandemic.
Particularly, we analyze the daily new cases of COVID-19 of representative contries for two stages: day $1$ (January 22${}^{\text{nd}}$ 2020) to day $100$ (April 30${}^{\text{th}}$ 2020) and day $101$ (May 1${}^{\text{st}}$ 2020) to day $391$ (February 15${}^{\text{th}}$ 2021).
Our continuity scaling is pairwisely applied to reconstruct the transmission causal network.
As shown in Fig.~\ref{E24mainpic},
China shows a significant effect on a few countries at the first stage
and this effect disappears at the second stage.
However, other countries show a different situation with China, whose external effect lasts as shown in Section III.E and Fig.~S18 of SI.
Our results accord with that China holds stringent epidemic control strategies with sporadic domestic infections,
as evidenced by official daily briefings,
demonstrating the potential of continuity scaling in detecting causal networks for ongoing complex systems.
Additionally,
We emphasize that day $100$ is a suitable critical day to distinguish the early severe stage and the late well-under-control stage of the pandemic [see Fig.~S18(a) of SI], while slight change of the critical day will not nullify our result.
As shown in Fig. S18(b) of SI,
when the critical day varies from day $94$ to day $106$,
no significant change (less than $5\%$) of the detected causal links occurs at both stages,
and the number of countries under influence of China at Stage $2$ remains zero.
See more details in Section III.E of SI.

Additional real world examples including air pollutants and hospital admission record from Hong Kong are also shown in Section III of SI.

\section{Discussion}

To summarize, we have developed a novel framework for data based detection and quantification of causation in complex dynamical systems.
On the basis of the widely used cross-map-based techniques, our framework enjoys a rigorous foundation, focusing on the continuity scaling law of the concerned system directly instead of only investigating the continuity of its cross map.
Therefore, our framework is consistent with the standard interpretation of causality, and works even in the typical cases where several existing typical methods do not perform that well or even they fail (see the comparison results in Section IV of SI).
In addition, the mathematical reasoning leading to the core of our framework, the continuity scaling, helps resolve the long-standing issue associated with techniques directly using cross map that information about the resulting variables is required to project the dynamical behavior of the causal variables,
whereas several works in the literature \cite{Quyen1999PhyD}, which directly studied the continuity or the smoothness of the cross map, likely yielded confused detected results on causal directions.

{\bf Computational complexity.} The computational complexity of the algorithm is $O(T^2N_\varepsilon)$, which is relatively smaller than the CCM method, whose computational complexity is $O(T^2\log T)$.

{\bf Limitations and future works.} Nevertheless, there are still some limitations in the presently proposed framework.   For instance, currently, only bivariate detection algorithm is designed, so generalization to multivariate network inference requires further considerations, as analogous to those works presented in Refs.~\cite{Peters2017elements,Runge2018causal,research2019}.  In addition, the causal time delay has not been taken into account in the current framework, so it also could be further investigated, similar to the work reported in Ref.~\cite{ma2017detection}.  Definitely, we will settle these questions in our future work.

Detecting causality in complex dynamical systems has broad applications not only in science and engineering, but also in many aspects of the modern society, demanding accurate, efficient, and rigorously justified and hence trustworthy methodologies.
Our present work provides a vehicle along this feat, and indeed resolves the puzzles arising in the use of those influential methods.

\section{Methods}
\footnotesize
{\bf Continuity scaling framework: a detailed description of algorithms.}
Let $\{u_t\}_{t=1,2,\dots,T}$ and $\{v_t\}_{t=1,2,\dots,T}$ be two
experimentally measured time series of internal variables $\{(\bm x_t,\bm y_t)\}_{t\in\mathbb{N}}$. Typically, if the dynamical variables $\{(\bm x_t,\bm y_t)\}_{t\in\mathbb{N}}$ are accessible, $\{(u_t,v_t)\}$ reduce to one-dimensional coordinate of the internal system. The key computational steps of our
continuity scaling framework are described, as follows.

\begin{shaded}
We reconstruct the phase space using the classical method of delay coordinate
embedding~\cite{Takens:1981} with the optimal embedding dimension $d_{z}$ and
time lag $\tau_{z}$ determined by the methods in
Refs.~\cite{MutualInfoMethod,FNNmethod} (i.e., the false nearest neighbors and the delayed mutual information respectively):
\begin{equation}
\mathcal{L}_z\triangleq\left\{\bm{z}(t)=\left(z_t,z_{t+\tau_z},\dots,z_{t+(d_z-1)\tau_z}\right)~|~ t=1,\dots,T_0\right\},
\end{equation}
where $z=u,v$, $T_0=\min\{T-(d_z-1)\tau_z~|~z=u,v\}$, and Euclidean distance
is used for both $\mathcal{L}_{u,v}$.
\end{shaded}

We present the steps for causation detection using the case of
$\bm{v}\hookrightarrow\bm{u}$ as an example.

\begin{shaded}
We calculate the respective diameters for $\mathcal{L}_{u,v}$ as
\begin{equation}\label{largest}
D_z \triangleq \max \left\{\text{dist}_{\mathcal{L}_{z}}(\bm{z}(t),\bm{z}(\tau))~\big|~1\leqslant t,\tau\leqslant T_0 \right\},
\end{equation}
where $z=u,v$, and $\bm{z}=\bm{u},\bm{v}$. We set up a group of numbers,
$\{\varepsilon_{\bm{u},j}\}_{j=1,\cdots,N_{\epsilon}}$, as
$\varepsilon_{\bm{u},1}=e\cdot D_u$, $\varepsilon_{\bm{u},N_{\epsilon}}=D_u$,
with the other elements satisfying
\begin{equation}
  (\ln \varepsilon_{\bm{u},j} -\ln \varepsilon_{\bm{u},1})/(j-1)=
(\ln \varepsilon_{\bm{u},N_{\varepsilon}} -\ln \varepsilon_{\bm{u},1})/(N_{\varepsilon}-1)
\end{equation}
for $j=2,\dots,N_{\varepsilon}-1$. Then, in light of \eqref{indexset} with \eqref{delta}, we have
\begin{equation}
\delta^t_{\bm{v}}(\varepsilon_{\bm{u}})= \#[I^t_{\bm{u}}(\varepsilon_{\bm{u}})]^{-1}
\sum_{\tau\in I^t_{\bm{u}}(\varepsilon_{\bm{u}})}
\text{dist}_{\mathcal{L}_{v}}(\bm{v}(t),\bm{v}(\tau-1)),
\end{equation}
with
\begin{equation}
\begin{array}{l}
I^t_{\bm{u}}(\varepsilon_{\bm{u}})=\{\tau\in\mathbb{N}~\big|~\text{dist}_{\mathcal{L}_{u}}(\bm{u}(t+1),\bm{u}(\tau))<\varepsilon_{\bm{u}},\\~~~~~~~~~~~~~~~~~~~~~~~~~~~~~~~~~~~~~~~~~~~~~~~~~|t+1-\tau|>E\},
\end{array}
\end{equation}
where, numerically, $\varepsilon_{\bm{u}}$ alters its value successively from
the set $\{\varepsilon_{\bm{u},j}\}_{j=1,\cdots,N_{\epsilon}}$, and the
threshold $E$ is a positive number chosen to avoid the situation where the
nearest neighboring points are induced by the consecutive time order only.
\end{shaded}

As defined, $\langle \delta^t_{\bm{v}}(\varepsilon_{\bm{u}}) \rangle_{t\in\mathbb{N}}$ is
the average of $\delta^t_{\bm{v}}(\varepsilon_{\bm{u}})$ over all possible
time $t$. We use a finite number of pairs
$\{(\langle\delta^{t}_{\bm{v}}(\varepsilon_{\bm{u},j})\rangle_{t\in\mathbb{N}_{T_{0}}},\ln\varepsilon_{\bm{u},j})\}_{j=1,\dots,N_{\varepsilon}}$
to approximate the scaling relation between
$\langle \delta^t_{\bm{v}}(\varepsilon_{\bm{u}})\rangle_{t\in\mathbb{N}}$
and $\ln \varepsilon_{\bm{u}}$, where $\mathbb{N}_{T_{0}}=\{1,2,\cdots,T_{0}\}$.
Theoretically, a larger value of $N_{\epsilon}$ and a smaller value of $e$
will result in a more accurate approximation of the scaling relation. In
practice, the accuracy is determined by the length of the observational time
series, the sampling duration, and different types of noise perturbations. In
numerical simulations, we set $e=0.001$ and $N_{\epsilon}=33$. In addition,
a too large or a too small value of $\varepsilon_{\bm u}$ can induce
insufficient data to restore the neighborhood and/or the entire manifold.
We thus set
$\delta^t_{\bm{v}}(\varepsilon_{\bm{u},j})=\delta^t_{\bm{v}}(\varepsilon_{\bm{u},j+1})$ as a practical technique as the number of points is limited practically in a small neighborhood. As a result, near zero slope values would
appear on both sides of the scaling curve
$\langle \delta^t_{\bm{v}}(\varepsilon_{\bm{u}})\rangle_{t\in\mathbb{N}}$-$\ln \varepsilon_{\bm u}$,
as demonstrated in Fig.~\ref{fig3} and in SI. In such a case, to estimate the
slope of the scaling relation, we take the following approach.

\begin{shaded}
Define a group of numbers by
\begin{equation}
S_{j}\triangleq
\frac{\langle\delta_{\bm{v}}^{t}(\varepsilon_{\bm{u},j+1})\rangle_{t\in\mathbb{N}_{T_{0}}}- \langle\delta^{t}_{\bm{v}}(\varepsilon_{\bm{u},j})\rangle_{t\in\mathbb{N}_{T_{0}}}}{\ln \varepsilon_{\bm{u},j+1}-\ln \varepsilon_{\bm{u},j}},
\end{equation}
where $j=1,\cdots,N_{\varepsilon}-1$, sort them in a descending order, from which we determine the
$[\frac{N_{\varepsilon}+1}{2}]$ largest numbers,
collect their subscripts - $j$'s together as an index set $\hat J$, and set
$H\triangleq\left\{j,~j+1~\big|~j\in \hat J \right\}$.
Applying the least squares method to the linear regression model:
\begin{equation}\label{lrq}
\langle\delta^{t}_{\bm{v}}(\varepsilon_{\bm{u}})\rangle_{t\in\mathbb{N}}=S \cdot \ln \varepsilon_{\bm{u}} +b
\end{equation}
with the dataset
$\{(\langle\delta^{t}_{\bm{v}}(\varepsilon_{\bm{u},j})\rangle_{t\in\mathbb{N}_{T_{0}}},\ln\varepsilon_{\bm{u},j})\}_{j\in H}$,
we get the optimal values $(\hat{S},\hat{b})$ for the parameters $(S,b)$
in \eqref{lrq} and finally obtain the slope of the scaling relation as
$s_{\bm{v} \hookrightarrow \bm{u}}\triangleq \hat{S}$.
\end{shaded}

For the other causal direction from $\bm{u}$ to $\bm{v}$, these steps are equally applicable to estimating the slope $s_{\bm{u}\hookrightarrow\bm{v}}$.

To assess the statistical significance of the numerically determined causation,
we devise the following surrogate test using the case of $\bm{v}$ causing
$\bm{u}$ as an illustrative example.

\begin{shaded}
Divide the time series $\{\bm{u}(t)\}_{t\in\mathbb{N}_{T_{0}}}$ into $N_G$
consecutive segments of equal length (except for the last segment - the
shortest segment). Randomly shuffle these segments and then regroup them into
a surrogate sequence $\{\bm{\hat{u}}(t)\}_{t\in\mathbb{N}_{T_{0}}}$. Applying
such a random permutation method to $\{\bm{v}(t)\}_{t\in\mathbb{N}_{T_{0}}}$
generates another surrogate sequence
$\{\bm{\hat{v}}(t)\}_{t\in\mathbb{N}_{T_{0}}}$. Carrying out the slope
computation yields $s_{\bm{\hat{v}} \hookrightarrow \bm{\hat{u}}}$. The
procedure can be repeated for a sufficient number of times, say $Q$, which
consequently yields a group of estimated slopes, denoted as
$\{s^{q}_{\bm{\hat{v}} \hookrightarrow \bm{\hat{u}}}\}_{q=0,1\cdots,Q}$,
where $s^{0}_{\bm{\hat{v}} \hookrightarrow \bm{\hat{u}}}$ is set as
$s_{\bm{v} \hookrightarrow \bm{u}}$ obtained from the original time series.
For all the estimated slopes, we calculate their mean $\hat{\mu}_{\bm{v}\hookrightarrow\bm{u}}$ and the standard deviation
$\hat{\sigma}_{\bm{v}\hookrightarrow\bm{u}}$. The $p$-value for
$s_{\bm{v} \hookrightarrow \bm{u}}$ is calculated as
\begin{equation}
p_{s_{\bm{v} \hookrightarrow \bm{u}}}\triangleq
1-\text{normcdf}\left[\frac{s_{\bm{v} \hookrightarrow \bm{u}}-
\hat{\mu}_{\bm{v}\hookrightarrow\bm{u}}}{\hat{\sigma}_{\bm{v}\hookrightarrow\bm{u}}}\right],
\end{equation}
where $\text{normcdf}[\cdot]$ is the cumulative Gaussian distribution function.
The principle of statistical hypothesis testing guarantees the existence of
causation from $\bm{v}$ to $\bm{u}$ if
$p_{s_{\bm{v} \hookrightarrow \bm{u}}} < 0.05$.
\end{shaded}

In simulations, we set the number of segments to be $N_G=25$ and the number
of times for random permutations to be $Q \ge 20$.

\section{Acknowledgements}

W.L. is supported by the National Key R\&D Program of China (Grant
No.~2018YFC0116600), by the National Natural Science Foundation of China (Grant Nos.~11925103 \& 61773125), by the STCSM (Grant No.~18DZ1201000), and by the Shanghai Municipal Science and Technology Major Project (No. 2021SHZDZX0103).
Y.-C.L. is supported by AFOSR (Grant No. FA9550-21-1-0438).
S.-Y.L. is supported by the National Natural Science Foundation of China (No.{\ }12101133) and ``Chenguang Program'' supported by Shanghai Education Development Foundation and Shanghai Municipal Education Commission (No.{\ }20CG01).
Q.N. is partially supported by
NSF (Grant No.~DMS1763272) and the Simons Foundation (Grant No.~594598).

\section{Competing Interests}

The authors declare no competing interests.

\section{Code Availability}
The source codes for our CS framework are available at \url{https://github.com/bianzhiyu/ContinuityScaling}.

\section{Contributions}
W.L. conceived idea;  X.Y., S.-Y.L., and W.L. designed and performed research;  X.Y., S.-Y.L., H.-F.M., and W.L. analyzed data,  H.-F.M., Y.-C.L., and Q.N. contributed data and analysis tools, and all the authors wrote the paper.   X.Y. and S.-Y.L. were equally contributed.


\begin{thebibliography}{0}%
\makeatletter
\providecommand \@ifxundefined [1]{%
 \@ifx{#1\undefined}
}%
\providecommand \@ifnum [1]{%
 \ifnum #1\expandafter \@firstoftwo
 \else \expandafter \@secondoftwo
 \fi
}%
\providecommand \@ifx [1]{%
 \ifx #1\expandafter \@firstoftwo
 \else \expandafter \@secondoftwo
 \fi
}%
\providecommand \natexlab [1]{#1}%
\providecommand \enquote  [1]{``#1''}%
\providecommand \bibnamefont  [1]{#1}%
\providecommand \bibfnamefont [1]{#1}%
\providecommand \citenamefont [1]{#1}%
\providecommand \href@noop [0]{\@secondoftwo}%
\providecommand \href [0]{\begingroup \@sanitize@url \@href}%
\providecommand \@href[1]{\@@startlink{#1}\@@href}%
\providecommand \@@href[1]{\endgroup#1\@@endlink}%
\providecommand \@sanitize@url [0]{\catcode `\\12\catcode `\$12\catcode
  `\&12\catcode `\#12\catcode `\^12\catcode `\_12\catcode `\%12\relax}%
\providecommand \@@startlink[1]{}%
\providecommand \@@endlink[0]{}%
\providecommand \url  [0]{\begingroup\@sanitize@url \@url }%
\providecommand \@url [1]{\endgroup\@href {#1}{\urlprefix }}%
\providecommand \urlprefix  [0]{URL }%
\providecommand \Eprint [0]{\href }%
\providecommand \doibase [0]{http://dx.doi.org/}%
\providecommand \selectlanguage [0]{\@gobble}%
\providecommand \bibinfo  [0]{\@secondoftwo}%
\providecommand \bibfield  [0]{\@secondoftwo}%
\providecommand \translation [1]{[#1]}%
\providecommand \BibitemOpen [0]{}%
\providecommand \bibitemStop [0]{}%
\providecommand \bibitemNoStop [0]{.\EOS\space}%
\providecommand \EOS [0]{\spacefactor3000\relax}%
\providecommand \BibitemShut  [1]{\csname bibitem#1\endcsname}%
\let\auto@bib@innerbib\@empty
\end{thebibliography}%


\begin{thebibliography}{10}
\expandafter\ifx\csname url\endcsname\relax
  \def\url#1{\texttt{#1}}\fi
\expandafter\ifx\csname urlprefix\endcsname\relax\def\urlprefix{URL }\fi
\providecommand{\bibinfo}[2]{#2}
\providecommand{\eprint}[2][]{\url{#2}}

\bibitem{bunge2017causality}
\bibinfo{author}{Bunge, M.}
\newblock \emph{\bibinfo{title}{Causality and modern science}}
  (\bibinfo{publisher}{Routledge}, \bibinfo{year}{2017}).

\bibitem{pearl2009causality}
\bibinfo{author}{Pearl, J.}
\newblock \emph{\bibinfo{title}{Causality}} (\bibinfo{publisher}{Cambridge
  university press}, \bibinfo{year}{2009}).

\bibitem{Runge2019inferring}
\bibinfo{author}{Runge, J.}, \bibinfo{author}{Bathiany, S.} \emph{et~al.}
\newblock \bibinfo{title}{Inferring causation from time series in Earth system sciences}.
\newblock \emph{\bibinfo{journal}{Nature communications}} \textbf{\bibinfo{volume}{10}}, \bibinfo{pages}{2553} (\bibinfo{year}{2019}).
\newblock doi: \href{https://www.doi.org/10.1038/s41467-019-10105-3}{10.1038/s41467-019-10105-3}.

\bibitem{CV:2015}
\bibinfo{author}{Collins, F.~S.} \& \bibinfo{author}{Varmus, H.}
\newblock \bibinfo{title}{A new initiative on precision medicine}.
\newblock \emph{\bibinfo{journal}{New England Journal of Medicine}}
  \textbf{\bibinfo{volume}{372}}, \bibinfo{pages}{793--795}
  (\bibinfo{year}{2015}).
\newblock doi: \href{https://www.doi.org/10.1056/NEJMp1500523}{10.1056/NEJMp1500523}.

\bibitem{SSFRLPWBBA:2016}
\bibinfo{author}{Saxe, G.~N.} \emph{et~al.}
\newblock \bibinfo{title}{A complex systems approach to causal discovery in
  psychiatry}.
\newblock \emph{\bibinfo{journal}{PloS ONE}} \textbf{\bibinfo{volume}{11}},
  \bibinfo{pages}{e0151174} (\bibinfo{year}{2016}).
\newblock doi: \href{https://www.doi.org/10.1371/journal.pone.0151174}{10.1371/journal.pone.0151174}.

\bibitem{CoxHinkley1979}
\bibinfo{author}{Cox, D.~R.} \& \bibinfo{author}{Hinkley, D.~V.}
\newblock \emph{\bibinfo{title}{Theoretical statistics}}
  (\bibinfo{publisher}{CRC Press}, \bibinfo{year}{1979}).
\newblock doi: \href{https://www.doi.org/10.1201/b14832}{10.1201/b14832}.

\bibitem{CoverThomas2012}
\bibinfo{author}{Cover, T.~M.}
\newblock \emph{\bibinfo{title}{Elements of information theory}}
  (\bibinfo{publisher}{John Wiley \& Sons}, \bibinfo{year}{1999}).


\bibitem{JPearl2009}
\bibinfo{author}{Pearl, J.} \emph{et~al.}
\newblock \bibinfo{title}{Causal inference in statistics: An overview}.
\newblock \emph{\bibinfo{journal}{Statistics Surveys}}
  \textbf{\bibinfo{volume}{3}}, \bibinfo{pages}{96--146}
  (\bibinfo{year}{2009}).
\newblock doi: \href{https://www.doi.org/10.1214/09-SS057}{10.1214/09-SS057}.

\bibitem{Wiener:1956}
\bibinfo{author}{Wiener, N.}
\newblock \bibinfo{title}{The theory of prediction}.
\newblock \emph{\bibinfo{journal}{Modern mathematics for engineers}}
  (\bibinfo{year}{1956}).

\bibitem{Granger:1969}
\bibinfo{author}{Granger, C.~W.}
\newblock \bibinfo{title}{Investigating causal relations by econometric models
  and cross-spectral methods}.
\newblock \emph{\bibinfo{journal}{Econometrica: Journal of the Econometric
  Society}} \bibinfo{pages}{424--438} (\bibinfo{year}{1969}).
\newblock doi: \href{https://www.doi.org/10.2307/1912791}{10.2307/1912791}.

\bibitem{Haufe2013}
\bibinfo{author}{Haufe, S.}{\ \textit{et al.}}
\newblock\bibinfo{title}{A critical assessment of connectivity measures for EEG data: A simulation study}.
\newblock\emph{\bibinfo{journal}{NeuroImage}}
\bibinfo{pages}{120--133}(\bibinfo{year}{2013}).
\newblock doi: \href{https://www.doi.org/10.1016/j.neuroimage.2012.09.036}{10.1016/j.neuroimage.2012.09.036}.


\bibitem{ding2006granger}
\bibinfo{author}{Ding, M.}, \bibinfo{author}{Chen, Y.} \&
  \bibinfo{author}{Bressler, S.~L.}
\newblock \bibinfo{title}{Granger causality: basic theory and application to
  neuroscience}.
\newblock \emph{\bibinfo{journal}{Handbook of time series analysis: recent
  theoretical developments and applications}} \textbf{\bibinfo{volume}{437}}
  (\bibinfo{year}{2006}).
\newblock doi: \href{https://www.doi.org/10.1002/9783527609970.ch17}{10.1002/9783527609970.ch17}.



\bibitem{Schreiber:2000}
\bibinfo{author}{Schreiber, T.}
\newblock \bibinfo{title}{Measuring information transfer}.
\newblock \emph{\bibinfo{journal}{Physical Review Letters}}
  \textbf{\bibinfo{volume}{85}}, \bibinfo{pages}{461} (\bibinfo{year}{2000}).
  \newblock doi: \href{https://www.doi.org/10.1103/PhysRevLett.85.461}{10.1103/PhysRevLett.85.461}.


\bibitem{frenzel2007partial}
\bibinfo{author}{Frenzel, S.} \& \bibinfo{author}{Pompe, B.}
\newblock \bibinfo{title}{Partial mutual information for coupling analysis of
  multivariate time series}.
\newblock \emph{\bibinfo{journal}{Physical Review Letters}}
  \textbf{\bibinfo{volume}{99}}, \bibinfo{pages}{204101}
  (\bibinfo{year}{2007}).
   \newblock doi: \href{https://www.doi.org/10.1103/PhysRevLett.99.204101}{10.1103/PhysRevLett.99.204101}.


\bibitem{vicente2011transfer}
\bibinfo{author}{Vicente, R.}, \bibinfo{author}{Wibral, M.},
  \bibinfo{author}{Lindner, M.} \& \bibinfo{author}{Pipa, G.}
\newblock \bibinfo{title}{Transfer entropy--a model-free measure of effective
  connectivity for the neurosciences}.
\newblock \emph{\bibinfo{journal}{Journal of Computational Neuroscience}}
  \textbf{\bibinfo{volume}{30}}, \bibinfo{pages}{45--67}
  (\bibinfo{year}{2011}).
   \newblock doi: \href{https://www.doi.org/10.1007/s10827-010-0262-3}{10.1007/s10827-010-0262-3}.


\bibitem{runge2012escaping}
\bibinfo{author}{Runge, J.}, \bibinfo{author}{Heitzig, J.},
  \bibinfo{author}{Petoukhov, V.} \& \bibinfo{author}{Kurths, J.}
\newblock \bibinfo{title}{Escaping the curse of dimensionality in estimating
  multivariate transfer entropy}.
\newblock \emph{\bibinfo{journal}{Physical Review Letters}}
  \textbf{\bibinfo{volume}{108}}, \bibinfo{pages}{258701}
  (\bibinfo{year}{2012}).
 \newblock doi: \href{https://www.doi.org/10.1103/PhysRevLett.108.258701}{10.1103/PhysRevLett.108.258701}.


\bibitem{SCB:2014}
\bibinfo{author}{Sun, J.}, \bibinfo{author}{Cafaro, C.} \&
  \bibinfo{author}{Bollt, E.~M.}
\newblock \bibinfo{title}{Identifying the coupling structure in complex systems
  through the optimal causation entropy principle}.
\newblock \emph{\bibinfo{journal}{Entropy}} \textbf{\bibinfo{volume}{16}},
  \bibinfo{pages}{3416--3433} (\bibinfo{year}{2014}).
   \newblock doi: \href{https://www.doi.org/10.3390/e16063416}{10.3390/e16063416}.



\bibitem{CSB:2015}
\bibinfo{author}{Cafaro, C.}, \bibinfo{author}{Lord, W.~M.},
  \bibinfo{author}{Sun, J.} \& \bibinfo{author}{Bollt, E.~M.}
\newblock \bibinfo{title}{Causation entropy from symbolic representations of
  dynamical systems}.
\newblock \emph{\bibinfo{journal}{Chaos: An Interdisciplinary Journal of
  Nonlinear Science}} \textbf{\bibinfo{volume}{25}}, \bibinfo{pages}{043106}
  (\bibinfo{year}{2015}).
 \newblock doi: \href{https://www.doi.org/10.1063/1.4916902}{10.1063/1.4916902}.




\bibitem{STB:2015}
\bibinfo{author}{Sun, J.}, \bibinfo{author}{Taylor, D.} \&
  \bibinfo{author}{Bollt, E.~M.}
\newblock \bibinfo{title}{Causal network inference by optimal causation
  entropy}.
\newblock \emph{\bibinfo{journal}{SIAM Journal on Applied Dynamical Systems}}
  \textbf{\bibinfo{volume}{14}}, \bibinfo{pages}{73--106}
  (\bibinfo{year}{2015}).
 \newblock doi: \href{https://www.doi.org/10.1137/140956166}{10.1137/140956166}.


\bibitem{research2021}
\bibinfo{author}{Solyanik-Gorgone, M.}, \bibinfo{author}{Ye, J.}, \bibinfo{author}{Miscuglio, M.}, \bibinfo{author}{Afanasev, A.}, \bibinfo{author}{Willner, A.~E.} \&
  \bibinfo{author}{Sorger, V.~J.}
\newblock \bibinfo{title}{Quantifying Information via Shannon Entropy in Spatially Structured Optical Beams}.
\newblock \emph{\bibinfo{journal}{Research}}
  \textbf{\bibinfo{volume}{2021}}, \bibinfo{pages}{9780760}
  (\bibinfo{year}{2021}).
 \newblock doi: \href{https://doi.org/10.34133/2021/9780760}{10.34133/2021/9780760}.


\bibitem{PRE_Hirata2010}
\bibinfo{autho}{Hirata, Y.}{\ }\&{\ }\bibinfo{author}{Aihara, K.}
\newblock \bibinfo{title}{Identifying hidden common causes from bivariate time series: A method using recurrence plots}.
\newblock\emph{\bibinfo{journal}{Physical Review E}}
 \textbf{\bibinfo{volume}{81}}, \bibinfo{pages}{016203} (\bibinfo{year}{2010}).
 \newblock doi: \href{https://www.doi.org/10.1103/PhysRevE.81.016203}{10.1103/PhysRevE.81.016203}.



\bibitem{PRE_Quiroga2000}
\bibinfo{author}{Quiroga, R.~Q.}, \bibinfo{author}{Arnhold, J.} \&
  \bibinfo{author}{Grassberger, P.}
\newblock \bibinfo{title}{Learning driver-response relationships from
  synchronization patterns}.
\newblock \emph{\bibinfo{journal}{Physical Review E}}
  \textbf{\bibinfo{volume}{61}}, \bibinfo{pages}{5142} (\bibinfo{year}{2000}).
 \newblock doi: \href{https://www.doi.org/10.1103/PhysRevE.61.5142}{10.1103/PhysRevE.61.5142}.



\bibitem{PhyD_Arnhold1999}
\bibinfo{author}{Arnhold, J.}, \bibinfo{author}{Grassberger, P.},
  \bibinfo{author}{Lehnertz, K.} \& \bibinfo{author}{Elger, C.~E.}
\newblock \bibinfo{title}{A robust method for detecting interdependences:
  application to intracranially recorded eeg}.
\newblock \emph{\bibinfo{journal}{Physica D: Nonlinear Phenomena}}
  \textbf{\bibinfo{volume}{134}}, \bibinfo{pages}{419--430}
  (\bibinfo{year}{1999}).
 \newblock doi: \href{https://www.doi.org/10.1016/S0167-2789(99)00140-2}{10.1016/S0167-2789(99)00140-2}.



\bibitem{harnack2017topological}
\bibinfo{author}{Harnack, D.}, \bibinfo{author}{Laminski, E.},
  \bibinfo{author}{Sch{\"u}nemann, M.} \& \bibinfo{author}{Pawelzik, K.~R.}
\newblock \bibinfo{title}{Topological causality in dynamical systems}.
\newblock \emph{\bibinfo{journal}{Physical Review Letters}}
  \textbf{\bibinfo{volume}{119}}, \bibinfo{pages}{098301}
  (\bibinfo{year}{2017}).
 \newblock doi: \href{https://www.doi.org/10.1103/PhysRevLett.119.098301}{10.1103/PhysRevLett.119.098301}.



\bibitem{SMYHDFM:2012}
\bibinfo{author}{Sugihara, G.} \emph{et~al.}
\newblock \bibinfo{title}{Detecting causality in complex ecosystems}.
\newblock \emph{\bibinfo{journal}{Science}} \textbf{\bibinfo{volume}{338}},
  \bibinfo{pages}{496--500} (\bibinfo{year}{2012}).
 \newblock doi: \href{https://www.doi.org/10.1126/science.1227079}{10.1126/science.1227079}.



\bibitem{DFHKMMPYS:2013}
\bibinfo{author}{Deyle, E.~R.} \emph{et~al.}
\newblock \bibinfo{title}{Predicting climate effects on pacific sardine}.
\newblock \emph{\bibinfo{journal}{Proceedings of the National Academy of
  Sciences}} \textbf{\bibinfo{volume}{110}}, \bibinfo{pages}{6430--6435}
  (\bibinfo{year}{2013}).
   \newblock doi: \href{https://www.doi.org/10.1073/pnas.1215506110}{10.1073/pnas.1215506110}.


\bibitem{WPCFMCHMPW:2014}
\bibinfo{author}{Wang, X.} \emph{et~al.}
\newblock \bibinfo{title}{A two-fold increase of carbon cycle sensitivity to
  tropical temperature variations}.
\newblock \emph{\bibinfo{journal}{Nature}} \textbf{\bibinfo{volume}{506}},
  \bibinfo{pages}{212--215} (\bibinfo{year}{2014}).
 \newblock doi: \href{https://www.doi.org/10.1038/nature12915}{10.1038/nature12915}.


\bibitem{MAC:2014}
\bibinfo{author}{Ma, H.}, \bibinfo{author}{Aihara, K.} \&
  \bibinfo{author}{Chen, L.}
\newblock \bibinfo{title}{Detecting causality from nonlinear dynamics with
  short-term time series}.
\newblock \emph{\bibinfo{journal}{Scientific Reports}}
  \textbf{\bibinfo{volume}{4}}, \bibinfo{pages}{1--10} (\bibinfo{year}{2014}).
 \newblock doi: \href{https://www.doi.org/10.1038/srep07464}{10.1038/srep07464}.


\bibitem{MW:2014}
\bibinfo{author}{McCracken, J.~M.} \& \bibinfo{author}{Weigel, R.~S.}
\newblock \bibinfo{title}{Convergent cross-mapping and pairwise asymmetric
  inference}.
\newblock \emph{\bibinfo{journal}{Physical Review E}}
  \textbf{\bibinfo{volume}{90}}, \bibinfo{pages}{062903}
  (\bibinfo{year}{2014}).
 \newblock doi: \href{https://www.doi.org/10.1103/PhysRevE.90.062903}{10.1103/PhysRevE.90.062903}.


\bibitem{YDGS:2015}
\bibinfo{author}{Ye, H.}, \bibinfo{author}{Deyle, E.~R.},
  \bibinfo{author}{Gilarranz, L.~J.} \& \bibinfo{author}{Sugihara, G.}
\newblock \bibinfo{title}{Distinguishing time-delayed causal interactions using
  convergent cross mapping}.
\newblock \emph{\bibinfo{journal}{Scientific Reports}}
  \textbf{\bibinfo{volume}{5}}, \bibinfo{pages}{14750} (\bibinfo{year}{2015}).
 \newblock doi: \href{https://www.doi.org/10.1038/srep14750}{10.1038/srep14750}.


\bibitem{clark2015spatial}
\bibinfo{author}{Clark, A.~T.} \emph{et~al.}
\newblock \bibinfo{title}{Spatial convergent cross mapping to detect causal
  relationships from short time series}.
\newblock \emph{\bibinfo{journal}{Ecology}} \textbf{\bibinfo{volume}{96}},
  \bibinfo{pages}{1174--1181} (\bibinfo{year}{2015}).
 \newblock doi: \href{https://www.doi.org/10.1890/14-1479.1}{10.1890/14-1479.1}.



\bibitem{JHHLL:2016}
\bibinfo{author}{Jiang, J.-J.}, \bibinfo{author}{Huang, Z.-G.},
  \bibinfo{author}{Huang, L.}, \bibinfo{author}{Liu, H.} \&
  \bibinfo{author}{Lai, Y.-C.}
\newblock \bibinfo{title}{Directed dynamical influence is more detectable with
  noise}.
\newblock \emph{\bibinfo{journal}{Scientific Reports}}
  \textbf{\bibinfo{volume}{6}}, \bibinfo{pages}{24088} (\bibinfo{year}{2016}).
 \newblock doi: \href{https://www.doi.org/10.1038/srep24088}{10.1038/srep24088}.


\bibitem{ma2017detection}
\bibinfo{author}{Ma, H.} \emph{et~al.}
\newblock \bibinfo{title}{Detection of time delays and directional interactions
  based on time series from complex dynamical systems}.
\newblock \emph{\bibinfo{journal}{Physical Review E}}
  \textbf{\bibinfo{volume}{96}}, \bibinfo{pages}{012221}
  (\bibinfo{year}{2017}).
 \newblock doi: \href{https://www.doi.org/10.1103/PhysRevE.96.012221}{10.1103/PhysRevE.96.012221}.


\bibitem{Amigo2018}
\bibinfo{author}{Amig{\'o}, J.~M.} \& \bibinfo{author}{Hirata, Y.}
\newblock \bibinfo{title}{Detecting directional couplings from multivariate
  flows by the joint distance distribution}.
\newblock \emph{\bibinfo{journal}{Chaos: An Interdisciplinary Journal of
  Nonlinear Science}} \textbf{\bibinfo{volume}{28}}, \bibinfo{pages}{075302}
  (\bibinfo{year}{2018}).
 \newblock doi: \href{https://www.doi.org/10.1063/1.5010779}{10.1063/1.5010779}.



\bibitem{wang2018detecting}
\bibinfo{author}{Wang, Y.} \emph{et~al.}
\newblock \bibinfo{title}{Detecting the causal effect of soil moisture on
  precipitation using convergent cross mapping}.
\newblock \emph{\bibinfo{journal}{Scientific Reports}}
  \textbf{\bibinfo{volume}{8}}, \bibinfo{pages}{1--8} (\bibinfo{year}{2018}).
 \newblock doi: \href{https://www.doi.org/10.1038/s41598-018-30669-2}{10.1038/s41598-018-30669-2}.

\bibitem{pcm}
\bibinfo{author}{Leng, S.} \emph{et~al.}
\newblock \bibinfo{title}{Partial cross mapping eliminates indirect causal  influences}.
\newblock \emph{\bibinfo{journal}{Nature Communications}}
  \textbf{\bibinfo{volume}{11}}, \bibinfo{pages}{1--9} (\bibinfo{year}{2020}).
 \newblock doi: \href{https://www.doi.org/10.1038/s41467-020-16238-0}{10.1038/s41467-020-16238-0}.



\bibitem{Takens:1981}
\bibinfo{author}{Takens, F.}
\newblock \bibinfo{title}{Detecting strange attractors in turbulence}.
\newblock In \emph{\bibinfo{booktitle}{Dynamical systems and turbulence,
  Warwick 1980}}, \bibinfo{pages}{366--381} (\bibinfo{publisher}{Springer},
  \bibinfo{year}{1981}).
 \newblock doi: \href{https://www.doi.org/10.1007/BFb0091924}{10.1007/BFb0091924}.



\bibitem{PCFS:1980}
\bibinfo{author}{Packard, N.~H.}, \bibinfo{author}{Crutchfield, J.~P.},
  \bibinfo{author}{Farmer, J.~D.} \& \bibinfo{author}{Shaw, R.~S.}
\newblock \bibinfo{title}{Geometry from a time series}.
\newblock \emph{\bibinfo{journal}{Physical Review Letters}}
  \textbf{\bibinfo{volume}{45}}, \bibinfo{pages}{712} (\bibinfo{year}{1980}).
 \newblock doi: \href{https://www.doi.org/10.1103/PhysRevLett.45.712}{10.1103/PhysRevLett.45.712}.



\bibitem{sauer1991embedology}
\bibinfo{author}{Sauer, T.}, \bibinfo{author}{Yorke, J.~A.} \&
  \bibinfo{author}{Casdagli, M.}
\newblock \bibinfo{title}{Embedology}.
\newblock \emph{\bibinfo{journal}{Journal of Statistical Physics}}
  \textbf{\bibinfo{volume}{65}}, \bibinfo{pages}{579--616}
  (\bibinfo{year}{1991}).
 \newblock doi: \href{https://www.doi.org/10.1007/BF01053745}{10.1007/BF01053745}.



\bibitem{stark1999delay}
\bibinfo{author}{Stark, J.}
\newblock \bibinfo{title}{Delay embeddings for forced systems. i. deterministic
  forcing}.
\newblock \emph{\bibinfo{journal}{Journal of Nonlinear Science}}
  \textbf{\bibinfo{volume}{9}}, \bibinfo{pages}{255--332}
  (\bibinfo{year}{1999}).
 \newblock doi: \href{https://www.doi.org/10.1007/s003329900072}{10.1007/s003329900072}.



\bibitem{stark2003delay}
\bibinfo{author}{Stark, J.}, \bibinfo{author}{Broomhead, D.~S.},
  \bibinfo{author}{Davies, M.~E.} \& \bibinfo{author}{Huke, J.}
\newblock \bibinfo{title}{Delay embeddings for forced systems. ii. stochastic
  forcing}.
\newblock \emph{\bibinfo{journal}{Journal of Nonlinear Science}}
  \textbf{\bibinfo{volume}{13}}, \bibinfo{pages}{519--577}
  (\bibinfo{year}{2003}).
 \newblock doi: \href{https://www.doi.org/10.1007/s00332-003-0534-4}{10.1007/s00332-003-0534-4}.


\bibitem{muldoon1998delay}
\bibinfo{author}{Muldoon, M.~R.}, \bibinfo{author}{Broomhead, D.~S.},
  \bibinfo{author}{Huke, J.~P.} \& \bibinfo{author}{Hegger, R.}
\newblock \bibinfo{title}{Delay embedding in the presence of dynamical noise}.
\newblock \emph{\bibinfo{journal}{Dynamics and Stability of Systems}}
  \textbf{\bibinfo{volume}{13}}, \bibinfo{pages}{175--186}
  (\bibinfo{year}{1998}).
 \newblock doi: \href{https://www.doi.org/10.1080/02681119808806259}{10.1080/02681119808806259}.


\bibitem{spirtes2000causation}
\bibinfo{author}{Spirtes, P.}, \bibinfo{author}{Glymour, C.~N.},
  \bibinfo{author}{Scheines, R.} \& \bibinfo{author}{Heckerman, D.}
\newblock \emph{\bibinfo{title}{Causation, prediction, and search}}
  (\bibinfo{publisher}{MIT press}, \bibinfo{year}{2000}).


\bibitem{Cummins2015}
\bibinfo{author}{Cummins, B.},
\bibinfo{author}{Gedeon, T.}{\ }\&
\bibinfo{author}{Spendlove, K.}
\newblock
\bibinfo{title}{On the efficacy of state space reconstruction methods in determining causality}.
\newblock \emph{\bibinfo{journal}{SIAM Journal on Applied Dynamical Systems}}
\textbf{\bibinfo{volume}{14}}, \bibinfo{pages}{335--381}
(\bibinfo{year}{2015}).
 \newblock doi: \href{https://www.doi.org/10.1137/130946344}{10.1137/130946344}.


\bibitem{kantz2004nonlinear}
\bibinfo{author}{Kantz, H.} \& \bibinfo{author}{Schreiber, T.}
\newblock \emph{\bibinfo{title}{Nonlinear time series analysis}},
  vol.~\bibinfo{volume}{7} (\bibinfo{publisher}{Cambridge university press},
  \bibinfo{year}{2004}).
 \newblock doi: \href{https://www.doi.org/10.1017/CBO9780511755798}{10.1017/CBO9780511755798}.


\bibitem{lancaster2018surrogate}
\bibinfo{author}{Lancaster, G.}, \bibinfo{author}{Iatsenko, D.},
  \bibinfo{author}{Pidde, A.}, \bibinfo{author}{Ticcinelli, V.} \&
  \bibinfo{author}{Stefanovska, A.}
\newblock \bibinfo{title}{Surrogate data for hypothesis testing of physical
  systems}.
\newblock \emph{\bibinfo{journal}{Physics Reports}}
  \textbf{\bibinfo{volume}{748}}, \bibinfo{pages}{1--60}
  (\bibinfo{year}{2018}).
 \newblock doi: \href{https://www.doi.org/10.1016/j.physrep.2018.06.001}{10.1016/j.physrep.2018.06.001}.



\bibitem{marbach2009generating}
\bibinfo{author}{Marbach, D.}, \bibinfo{author}{Schaffter, T.},
  \bibinfo{author}{Mattiussi, C.} \& \bibinfo{author}{Floreano, D.}
\newblock \bibinfo{title}{Generating realistic in silico gene networks for
  performance assessment of reverse engineering methods}.
\newblock \emph{\bibinfo{journal}{Journal of Computational Biology}}
  \textbf{\bibinfo{volume}{16}}, \bibinfo{pages}{229--239}
  (\bibinfo{year}{2009}).
 \newblock doi: \href{https://www.doi.org/10.1089/cmb.2008.09TT}{10.1089/cmb.2008.09TT}.



\bibitem{marbach2010revealing}
\bibinfo{author}{Marbach, D.} \emph{et~al.}
\newblock \bibinfo{title}{Revealing strengths and weaknesses of methods for
  gene network inference}.
\newblock \emph{\bibinfo{journal}{Proceedings of the National Academy of
  Sciences}} \textbf{\bibinfo{volume}{107}}, \bibinfo{pages}{6286--6291}
  (\bibinfo{year}{2010}).
 \newblock doi: \href{https://www.doi.org/10.1073/pnas.0913357107}{10.1073/pnas.0913357107}.


\bibitem{Lasker1983Ocean}
\bibinfo{author}{R. Lasker} \& \bibinfo{author}{A. MacCall}.
\newblock \bibinfo{title}{New ideas on the fluctuations of the clupeoid stocks off California}.
\newblock In \bibinfo{booktitle}{\emph{Proceedings of the Joint Oceanographic Assembly, Halifax, August 1982: General Symposia}}, \bibinfo{pages}{110--120} (\bibinfo{publisher}{Department of Fisheries and Oceans, Ontario, 1983}).


\bibitem{Quyen1999PhyD}
\bibinfo{author}{Quyen, M.~L.~V.}, \bibinfo{author}{Martinerie, J.},
\bibinfo{author}{Adam, C.}, \bibinfo{author}{Varela, F.~J.}
\newblock \bibinfo{title}{Nonlinear analyses of interictal EEG map the brain interdependences in human focal epilepsy}.
\newblock \emph{\bibinfo{journal}{Physica D: Nonlinear Phenomena}}
\textbf{\bibinfo{volume}{127}},\bibinfo{pages}{250--266}(\bibinfo{year}{1999}).
 \newblock doi: \href{https://www.doi.org/10.1016/S0167-2789(98)00258-9}{10.1016/S0167-2789(98)00258-9}.


\bibitem{Peters2017elements}
\bibinfo{author}{Peters, J.}, \bibinfo{author}{Janzing, D.} \emph{et~al.}
\newblock \emph{\bibinfo{title}{Elements of causal inference: foundations and learning algorithms}}
  (\bibinfo{publisher}{MIT Press}, \bibinfo{year}{2017}).


\bibitem{Runge2018causal}
\bibinfo{author}{Runge, J.}
\newblock \bibinfo{title}{Causal network reconstruction from time series: From theoretical assumptions to practical estimation}.
\newblock \emph{\bibinfo{journal}{Chaos: An Interdisciplinary Journal of Nonlinear Science}} \textbf{\bibinfo{volume}{28}}, \bibinfo{pages}{075310} (\bibinfo{year}{2018}).
 \newblock doi: \href{https://www.doi.org/10.1063/1.5025050}{10.1063/1.5025050}.


\bibitem{research2019}
\bibinfo{author}{Lou, Y.}, \bibinfo{author}{Wang, L.} \& \bibinfo{author}{Chen, G.}
\newblock \bibinfo{title}{Enhancing Controllability Robustness of q-Snapback Networks through Redirecting Edges}.
\newblock \emph{\bibinfo{journal}{Research}}
  \textbf{\bibinfo{volume}{2019}}, \bibinfo{pages}{7857534} (\bibinfo{year}{2019}).
 \newblock doi: \href{https://doi.org/10.34133/2019/7857534}{10.34133/2019/7857534}.


\bibitem{MutualInfoMethod}
\bibinfo{author}{Fraser, A.~M.} \& \bibinfo{author}{Swinney, H.~L.}
\newblock \bibinfo{title}{Independent coordinates for strange attractors from
  mutual information}.
\newblock \emph{\bibinfo{journal}{Physical Review A}}
  \textbf{\bibinfo{volume}{33}}, \bibinfo{pages}{1134} (\bibinfo{year}{1986}).
 \newblock doi: \href{https://www.doi.org/10.1103/PhysRevA.33.1134}{10.1103/PhysRevA.33.1134}.


\bibitem{FNNmethod}
\bibinfo{author}{Kennel, M.~B.}, \bibinfo{author}{Brown, R.} \&
  \bibinfo{author}{Abarbanel, H.~D.}
\newblock \bibinfo{title}{Determining embedding dimension for phase-space
  reconstruction using a geometrical construction}.
\newblock \emph{\bibinfo{journal}{Physical Review A}}
  \textbf{\bibinfo{volume}{45}}, \bibinfo{pages}{3403} (\bibinfo{year}{1992}).
 \newblock doi: \href{https://www.doi.org/10.1103/PhysRevA.45.3403}{10.1103/PhysRevA.45.3403}.


\end{thebibliography}
\end{document}